\documentclass[11pt,leqno]{article}
 \usepackage[utf8,latin1]{inputenc}
 \usepackage[T1]{fontenc}
   \usepackage[english]{babel}
    \usepackage{lmodern}
 \usepackage{amsmath, amssymb, amsthm}

\def\R{{\mathbb R}}

\def\E{{\mathbb E}}

\newtheorem{theoreme}{Theorem}

\newtheorem{definition}{Definition}
\newtheorem{corollaire}{Corollary}

\begin{document}

\begin{center}
{\bf{\Large On Wintgen ideal submanifolds satisfying some
pseudo-symmetry type curvature conditions}}
\end{center}


\begin{center}
Ryszard Deszcz, Ma\l gorzata G\l ogowska, Miroslava Petrovi\'{c}-Torga\v{s}ev
and Georges Zafindratafa
\end{center}


\begin{center}
{\sl{Dedicated to Professor Bang-Yen Chen on the occasion of his 80th birthday}} 
\end{center}


\vspace{3mm}

\noindent
{\bf{Abstract.}}
Let $M$ be a Wintgen ideal submanifold of dimension $ n $
in a real space form $\R^{n+m}({\widetilde k})$ of dimension $(n+m)$
and of constant curvature $\widetilde k$, $n \geq 4$, $m \geq 1$.
Let $g$, $R$, $\rm Ricc$, $g \wedge \rm Ricc$ and $C$ be the metric tensor, 
the Riemann-Christoffel curvature tensor,
the Ricci tensor, the Kulkarni-Nomizu product of $g$ and $\rm Ricc$,
and the Weyl conformal curvature tensor of $M$, respectively.
In this paper we study Wintgen ideal submanifolds $M$ 
in real space forms $\R^{n+m}({\widetilde k})$, 
$n \geq 4$, $m \geq 1$,
satisfying  the following pseudo-symmetry type curvature conditions: 
\newline
(i) \ the tensors $ R \cdot C  $ and $ Q(g , R) $
(resp., $ Q(g , C) $, $ Q(g , g\wedge\mbox{\rm Ricc}) $,
$ Q(\mbox{\rm Ricc} , R) $ or $ Q(\mbox{\rm Ricc} ,
g\wedge\mbox{\rm Ricc}) $)  are linearly dependent;
\newline
(ii)  \ the tensors  $ C \cdot R  $ and $ Q(g , R) $
(resp., $ Q(g , C) $, $ Q(g , g\wedge\mbox{\rm Ricc}) $,
$ Q(\mbox{\rm Ricc} , R) $ or $ Q(\mbox{\rm Ricc} ,
g\wedge\mbox{\rm Ricc}) $)  are linearly dependent;
\newline
(iii) \  the tensors $ R \cdot C - C \cdot R  $ and $ Q(g , R) $
(resp., $ Q(g , C) $, $ Q(g , g\wedge\mbox{\rm Ricc}) $,
$ Q(\mbox{\rm Ricc} , R) $ or $ Q(\mbox{\rm Ricc} ,
g\wedge\mbox{\rm Ricc}) $)  are linearly dependent.


\bigskip

\tableofcontents


 \newpage

\section{Preliminaries}

\subsection{Pseudo-symmetry type curvature conditions}

Let $ M = M^n $ be a connected $n$-dimensional, $ n \geq 4 $,
Riemannian manifold of class $ C^{\infty}$.
Further, let $\nabla$ be its Levi-Civita connection
and $\mathfrak{X} (M)$ the Lie algebra of vector fields on $M$.
We define on $M$ the endomorphisms
$X \wedge _{A} Y$ and ${\mathcal{R}}(X,Y)$ of $\mathfrak{X} (M)$ by
$$
\begin{aligned}
(X \wedge _{A} Y)Z & = A(Y,Z)X - A(X,Z)Y, \\
{\mathcal R}(X,Y)Z & = \nabla _X \nabla _Y Z 
- \nabla _Y \nabla _X Z - \nabla _{[X,Y]}Z ,\\
\end{aligned}
$$
respectively,
where $X, Y, Z \in \mathfrak{X} (M)$
and
$A$ is a symmetric $(0,2)$-tensor on $M$.
The Ricci tensor $\mbox{\rm Ricc}$, the Ricci operator ${\mathcal{S}}$,
the scalar curvature $\kappa $ and
the endomorphism ${\mathcal{C}}(X,Y)$ of $(M,g)$
are defined by
$$
\begin{aligned}
\mbox{\rm Ricc}(X,Y)   & = \mathrm{tr} \{ Z \rightarrow {\mathcal{R}}(Z,X)Y \} ,\\
g({{\mathcal S}U} , V) & =  \mbox{\rm Ricc}(U , V)
     \  = \ \sum_{i=1}^n g({\mathcal R}(E_i,U)V , E_i) ,\\
\kappa & =  \mathrm{tr}\, {\mathcal{S}} ,\\
{\mathcal C}(X,Y)Z  & = {\mathcal R}(X,Y)Z
- \frac{1}{n-2}(X \wedge _{g} {\mathcal S}Y + {\mathcal S}X \wedge _{g} Y
- \frac{\kappa}{n-1}X \wedge _{g} Y)Z ,\\
\end{aligned}
$$
respectively,
where $X, Y, Z, U, V \in \mathfrak{X} (M)$ and
$\displaystyle \left\{ E_k \right\}_{k \in \left\{1, \ldots , n \right\}} $
is a local orthonormal tangent framefield.
The  {\sl  Riemann-Christoffel tensor} $R$ and
the {\sl  Weyl conformal curvature tensor} $C$ of $ M$ are defined by
$$
\begin{aligned}
	R(X  , Y, Z, W) &= g\left({\mathcal R}(X  , Y)Z , W \right), \\
	C(X  , Y, Z, W) &= g\left({\mathcal C}(X  , Y)Z , W \right), \\	
\end{aligned}
$$
respectively, where $X, Y, Z, W \in \mathfrak{X} (M)$.
Let $A$ and $B$ be the $(0,2)$-tensors;
their  {\sl  Kulkarni-Nomizu product} $A\wedge B$ is defined by
$$
\begin{aligned}
	(A\wedge B)(X_1 , X_2~; X , Y)
&= A(X_1 , Y) B(X_2 , X) + A(X_2 , X) B(X_1 , Y) \\
& \hspace{0,5cm}	- A(X_1 , X) B(X_2 , Y) - A(X_2 , Y) B(X_1 , X)  \\
\end{aligned}
 $$
where $ X, Y, X_1, X_2 \in \mathfrak{X} (M)$.
Now the tensor $C$ can be expressed in the form
$$
C = R - \frac{1}{n-2}{g\wedge\mbox{\rm Ricc}}
+ \frac{\kappa}{2(n-1)(n-2)}{g\wedge g} .
$$
Let $A$ be a symmetric $(0,2)$-tensor and $T$
a $(0,k)$-tensor, $ k \geq 1$~; we define
the {\sl  $(0,k+2)$-tensors}  $R \cdot T$ and  $Q(A , T)$ by
$$
\begin{aligned}
(R \cdot T)(X_1 , X_2 ,\dots , X_k ~; X , Y)
&= (R(X,Y) \cdot T)(X_1 , \dots , X_k) \\
&=	- T\left({\mathcal R}(X,Y)X_1 , X_2 ,\dots , X_k\right) \\
	& \hspace{0,5cm} - \cdots -  T\left(X_1 , X_2 ,\dots , X_{k-1} ,
	 {\mathcal R}(X,Y)X_k\right) \\
Q(A \cdot T)(X_1 , X_2 ,\dots , X_k ~; X , Y)
&= ((X\wedge_A Y) \cdot T)(X_1 , \dots , X_k) \\
&=	- T\left((X\wedge_A Y)X_1 , X_2 ,\dots , X_k\right) \\
	& \hspace{0,5cm} - \cdots -  T\left(X_1 , X_2 ,\dots , X_{k-1} ,
	  (X\wedge_A Y)X_k\right) \\
\end{aligned}
$$
respectively,
where $ X, Y, X_1, X_2, \ldots , X_n \in \mathfrak{X} (M)$.
If we set in the above formulas
$T=R$, $T=  \mbox{\rm Ricc} $, $T=C$, $A=g$ or $A=\mbox{\rm  Ricc}$,
 then we obtain  the tensors
$R \cdot R$, $ R \cdot \mbox{\rm Ricc}$,  $ R\cdot C$,
$ Q(g, R)$, $ Q(g, \mbox{\rm Ricc})$, $ Q(g, C)$,
 $ Q(g, g\wedge\mbox{\rm Ricc})$, $ Q(\mbox{\rm Ricc}, R)$
 and $ Q(\mbox{\rm Ricc}, C)$.
 Using the above definitions,
  we also define the tensors
 $   C\cdot R$, $ C\cdot C$   and $ C \cdot \mbox{\rm Ricc}$.
Further we have (see, e.g. \cite{{DGHS2011}, {DP-TVZ}})
$$
(n-2)({R\cdot C} - {C\cdot R}) =
Q\left(\mbox{\rm Ricc} - \frac{\kappa}{n-1}g ,  R\right)
- g\wedge (R \cdot \mbox{\rm Ricc}) + P
 $$
where the $(0,6)$-tensor $P$  is defined by
$$
\begin{aligned}
P(X_1, X_2, X_3, X_4 ~; X , Y)
&= g(X,X_1)R({\mathcal S}(Y), X_2 , X_3, X_4)
     - g(Y,X_1)R({\mathcal S}(X), X_2 , X_3, X_4)  \\
&\hspace{0,5cm} + g(X,X_2)R(X_1, {\mathcal S}(Y),   X_3, X_4)
     - g(Y,X_2)R(X_1, {\mathcal S}(X),  , X_3, X_4)  \\
&\hspace{0,5cm} + g(X,X_3)R(X_1, X_2,  {\mathcal S}(Y), X_4)
     - g(Y,X_3)R(X_1, X_2,  {\mathcal S}(Y), X_4)  \\
&\hspace{0,5cm} + g(X,X_4)R(X_1, X_2, X_3,  {\mathcal S}(Y))
     - g(Y,X_4)R(X_1, X_2, X_3, {\mathcal S}(X))  ,\\
\end{aligned}
$$
where $ X, Y, X_1, X_2,  X_3, X_4 \in \mathfrak{X} (M)$.

A Riemannian manifold $ (M , g)$ is said to be
{\sl semi-symmetric} (see \cite{Sz1982}) if
$R \cdot R = 0$ on $M$.
An extension form for semi-symmetric manifolds
are the pseudo-symmetric manifolds.
A Riemannian manifold $ (M , g)$ is said to be
{\sl pseudo-symmetric} (see \cite{{DES2}, {DGHS2011}}) if
at every point of $M$ the tensors
$R \cdot R$ and $Q(g,R)$ are linearly dependent,
which means that there exists a function
$L_R$ on the set
$ U_R = \left\{x \in M \ \big\vert \ R
- ( \kappa  / (2 n(n-1))) g \wedge g \ne 0  \ \text{at} \  x\right\} $
such that
\begin{equation}\label{pseudoSym}
R \cdot R = L_R Q(g , R) 
\end{equation}
on $U_R$.
Every semi-symmetric manifold is pseudo-symmetric. 
The converse is not true.
The warped product $S^1 \times_F S^{n-1} $ with
a positive smooth function $F$ is also pseudo-symmetric.
Pseudo-symmetric manifolds also are called
{\sl{pseudo-symmetric (in the sense of Deszcz)}}
or {\sl{Deszcz symmetric spaces}}
(see, e.g., \cite{{BYCH2021}, {15}, {LV1}, {LV2}, {LV3-Foreword}, {LV4}}).

A Riemannian manifold $ (M , g)$ is said to be
{\sl Weyl-semi-symmetric} (see \cite{{DES2}, {DGHS2011}}) if
$$
R \cdot C = 0 
$$
on $M$.
An extension form for
Weyl-semi-symmetric manifolds
are the Weyl pseudo-symmetric manifolds.
 A Riemannian manifold $ (M ,  g)$ is said to be
{\sl Weyl-pseudo-symmetric} (see  \cite{{DES2}, {DGHS2011}}) if
at every point of $M$ the tensors $R\cdot C$
and $Q(g,C)$ are linearly dependent,
which means that there exists a function
$ L_C$ on the set
$ U_C = \left\{x \in M \ \big\vert \   C  \ne 0 \ \text{at} \  x\right\} $
such that
\begin{equation}\label{WeylpseudoS}
R \cdot C = L_C Q(g , C) 
\end{equation}
on $U_C$.
It is obvious that every pseudo-symmetric manifold
is Weyl-pseudo-symmetric. The converse is not true.

A Riemannian manifold $ (M , g)$ is said to be
a {\sl manifold with  pseudo-symmetric Weyl tensor}
 (see   \cite{{DES2}, {DGHHY}, {DGHS2011}}) if
at every point of $M$ the tensors
$C \cdot C$ and $Q(g,C)$ are linearly dependent,
which means that there exists a function
$ L$ on the set
$U_C$
such that
\begin{equation}\label{pseudoWeylT}
C \cdot C = L Q(g , C) 
\end{equation}
on $U_C$.
Every Chen ideal submanifold in a space of constant curvature
is a  manifold with  pseudo-symmetric Weyl tensor
\cite{ {DGPV},  {2008_DP-TVZ}, {DP-TVZ}, {DeVerYap}}.

A Riemannian manifold $ (M , g)$ is said to be
{\sl Ricci-pseudo-symmetric}
(see  \cite{{DES2}, {DGHS2011}}) if  at every point of $M$
the tensors $R \cdot \mbox{\rm Ricc}$
and $Q(g,\mbox{\rm Ricc})$ are linearly dependent,
which means that there exists a function
$ L_{\mbox{\rm Ricc}}$ on the set
$ U_{\mbox{\rm Ricc}} =
\left\{x \in M \ \big\vert \ \mbox{\rm Ricc}
- \frac{\kappa}{n}g  \ne 0  \ \text{at} \ x\right\} $
such that
\begin{equation}\label{Riccipseudo}
R \cdot {\mbox{\rm Ricc}} = L_{\mbox{\rm Ricc}} Q(g , \mbox{\rm Ricc}) 
\end{equation}
on $ U_{\mbox{\rm Ricc}}$.
Every Cartan hypersurface $M$  in the sphere $ S^{n+1}$,
$ n=6,\ 12 \ \text{or} \ 24$,
is a non-pseudo-symmetric, Ricci-pseudo-symmetric with
non-pseudo-symmetric Weyl tensor \cite{DY} (see also \cite{44}).
We mention that
Ricci-pseudo-symmetric manifolds also are called
{\sl{Ricci pseudo-symmetric in the sense of Deszcz}},
or simply {\sl{Deszcz Ricci-symmetric}} (see, e.g., \cite{BYCH2021}).

A Riemannian manifold $ (M , g)$ is said to be
{\sl Ricci-Weyl-pseudo-symmetric} (see  \cite{{DES2}, {DGHS2011}})
 if at every point of $M$
the tensors $R \cdot C$
and $Q(\mbox{\rm Ricc},C)$ are linearly dependent,
which means that there exists a function
$ L'$ on the set
$ U  = \left\{x \in M  \ \big\vert \
Q(\mbox{\rm Ricc},C)  \ne 0 \ \text{at} \  x\right\} $
such that
   \begin{equation}\label{RicciWeyl}
	R \cdot C = L' Q(\mbox{\rm Ricc} , C) 
   \end{equation}
on $U$.
It is obvious that every semi-symmetric manifold ($R \cdot R = 0$)
satisfies (\ref{RicciWeyl}) with $L' = 0$.

We refer to \cite{{2023_DGHP-TZ}} (see also \cite{{2023_b_DGHP-TZ}, {DGHS-2022}})
for a recent survey on manifolds satisfying
(\ref{pseudoSym})--(\ref{RicciWeyl}) and other conditions
of this kind. Such conditions are called pseudo-symmetry type curvature conditions.
It seems that (\ref{pseudoSym}) is the most important condition
of that family of curvature conditions
(see \cite{ {DGHHY}, {DGP-TV02},  {DGPSS}, {DESZ}, {DP-TVZ},
{HaVer01}, {HaVer02}, {LV1}, {V2}, {LV2}, {LV3-Foreword}, {LV4}}).

\subsection{Submanifolds in space forms}

Throughout all sections of the present article,
 let $ M = M^n $   be a connected
 Riemannian manifold of class $ C^{\infty}$
	of  dimension $ n $
	in a  real space form   	$   \R^{n+m}({\widetilde k}) $
 of 	dimension $ (n+m) $
and of constant curvature $ {\widetilde k} $,
 $ n \geq 4 $, $ m \geq 1 $.

On $ \R^{n+m}({\widetilde k})  $, we denote by
  $ \widetilde{g} $ and   $ \widetilde{\nabla} $
respectively	 the {\sl Riemannian metric} and		
	 the cor\-res\-pon\-ding
		{\sl Levi-Civita connection}.
On the submanifold $ M $, the {\sl induced Riemannian metric}
and the corresponding  {\sl Levi-Civita connection}
 on $ M $ will be denoted by  $ g $,  $  \nabla $.
We will write as $ X $, $ Y $, $ \ldots  $
the {\sl tangent vector fields} on $ M $,
and as $ \xi $, $ \eta $, $ \ldots  $
the {\sl normal vector fields} on $ M $.

The well-known
{\sl formulae of Gauss}
and {\sl Weingarten} are given by
$$
\begin{aligned}
{\widetilde \nabla}_XY   & = \nabla_XY + h(X, Y),\\
{\widetilde \nabla}_X\xi & = \nabla^\perp_X\xi - A_\xi Y ,\\
\end{aligned}
$$
respectively,
whereby
$ \nabla^\perp $ is the {\sl normal connection}
induced in the
{\sl normal bundle} of $ M $ in $ \R^{n+m}({\widetilde k})  $,
$ h $ is the {\sl second fundamental form} of the submanifold $ M $
and
$ A_\xi  $ is the {\sl shape operator} or the
{\sl Weingarten map} on $ M $ with respect to the
normal vector field $ \xi $. We have
$$
 \displaystyle g(A_\xi(X), Y)  =  g\left(\xi, h(X,Y)\right)    ,$$
or still 
$$  \displaystyle  h(X,Y) = \sum_{\alpha=1}^m  g(A_\alpha(X), Y)\xi_\alpha  ,	$$
whereby
$\displaystyle \left\{\xi_\alpha\right\}_{\alpha \in \left\{1,
\ldots , m\right\}} $
is any {\sl local orthonormal normal framefield} on $ M $
in $ \R^{n+m}({\widetilde k})  $ and
whereby we have put
$ \displaystyle  A_\alpha =   A_{\xi_\alpha} $.

The {\sl mean curvature vector field} of  $ M $
in $ \R^{n+m}({\widetilde k})   $ is defined by
$$
\stackrel{\rightarrow}{H}
  =  \dfrac{1}{n}	\sum_{i=1}^n h(E_i,E_i)
= \dfrac{1}{n}\sum_{\alpha=1}^m \left(\mbox{\rm tr }A_{\alpha}\right)\xi_\alpha ,
$$
whereby
$  \displaystyle \left\{E_s\right\}_{s \in \left\{1, \ldots , n\right\}} $
is any {\sl local orthonormal tangent framefield} on $ M $, and
$  \displaystyle \left\{\xi_\alpha\right\}_{\alpha \in \left\{1,
\ldots , m\right\}} $
is any {\sl local orthonormal normal framefield} on $ M $
in $ \R^{n+m}({\widetilde k}) $~;
and its length $ H = \left\|\vec{H}\right\|$
is the {\sl mean curvature} of $M$ in $ \R^{n+m}({\widetilde k})  $.

The submanifold $ M $ in $ \R^{n+m}({\widetilde k}) $
is {\sl totally geodesic} when $ h = 0 $.
 $ M $ is {\sl totally umbilical} when $ h = g\stackrel{\rightarrow}{H} $.
It is {\sl minimal} when $ \stackrel{\rightarrow}{H}
= \stackrel{\rightarrow}{0} $,
	or equivalently, when
	its {\sl squared mean curvature function}
$H^2 =   \widetilde{g}\left(\stackrel{\rightarrow}{H} ,
\stackrel{\rightarrow}{H}\right)  $
	vanishes identically.
	$ M $ is   {\sl pseudo-umbilical} when the mean curvature vector field
	$  \stackrel{\rightarrow}{H} $ determines an
	{\sl umbilical normal direction} on $ M $
in $ \R^{n+m}({\widetilde k})  $, i.e.  when
$ A_{\stackrel{\rightarrow}{H}} = \chi \mbox{Id} $,
whereby   $ \mbox{Id} $ stands for the identity operator
on $ TM $ and $ \chi $ is some real function on $ M $.

Let
 $ R $ denote the {\sl induced  Riemann-Christoffel
 curvature}  corresponding  to the induced Levi-Civita connection
$  \nabla $ on $ M  $. Then according to the
{\sl equation of Gauss}
$$
	R(X,  Y, Z, W)
	=  \widetilde{g} \left(h(X, W)h(Y, Z) -  h(X, Z)h(Y, W)\right) 
	+  {\widetilde k}\left(g(X, W)g(Y, Z) -  g(X, Z)g(Y, W)\right)  ,
$$
where $X, Y, Z, W$ are tangent vector fields in $M$.

Let $\kappa$ be the {\sl scalar curvature function}
of $M$ in $ \R^{n+m}({\widetilde k})  $~; we have
$$
\kappa(p) = \sum_{i<j} K\left(p, E_i(p)\wedge E_j(p)\right)
$$
where $ K\left(p, E_i(p)\wedge E_j(p)\right)$
is the sectional curvature of $M$ at
a point $p \in M$ for the  plane section
$ \phi = E_i(p)\wedge E_j(p)$ in $ T_p M$
(where $\left\{E_i(p) ,  E_j(p)\right\}$ is independent).
By $  \inf(K)$ we will further denote the
{\sl function}
$$
\inf(K) :   M \to \R
$$
which attains to $p \in M$ 	
 the minimal value  $ \inf(K)(p)$ of all sectional
curvatures of $M $ at $p$.
The
{\sl normalized scalar curvature} of  $ (M , g) $ is defined as
follows
  $$
 \rho
  = \displaystyle   \dfrac{2}{n(n-1)}\sum_{i < j}
  R\left(E_i , E_j , E_j , E_i\right) ,
$$
whereby
$  \displaystyle \left\{E_s\right\}_{s \in \left\{1, \ldots , n\right\}} $
is any {\sl local orthonormal tangent framefield} on $ M $.

By the {\sl equation of Ricci},
the {\sl normal curvature tensor} $ R^\perp $
of $ M $ in  $ \R^{n+m}({\widetilde k})  $, i.e.,
the {\sl curvature tensor of the normal connection}
$ {\nabla}^\perp $ of $ M $ in   $ \R^{n+m}({\widetilde k}) $,
 is defined as follows
\begin{equation} \label{normcurvtens}
R^\perp (X , Y ; \xi , \eta)   =
\widetilde{g}\left({\nabla}^\perp_X{\nabla}^\perp_Y\xi -
 {\nabla}^\perp_Y{\nabla}^\perp_X\xi
- {\nabla}^\perp_{\left[X,Y\right]}\xi , \eta\right)
 = g\left(\left[A_\xi , A_\eta\right]X , Y\right) ,
\end{equation}
whereby
$ \left[A_\xi , A_\eta\right] =  A_\xi A_\eta  - A_\eta A_\xi $.

The
{\sl normalized scalar normal  curvature} $  \rho^\perp $  of  $ (M , g) $
    in    $ \R^{n+m}({\widetilde k})  $
  is given as follows
  $$
 \rho^\perp   = \displaystyle
  \dfrac{2}{n(n-1)}\, \sqrt{\sum_{i < j}\sum_{\alpha < \beta}
   R^\perp\left(E_i , E_j , \xi_\alpha , \xi_\beta\right)} ,
$$
whereby
$  \displaystyle \left\{E_s\right\}_{s \in \left\{1, \ldots , n\right\}} $
is any {\sl local orthonormal tangent framefield} on $ M $, and
$  \displaystyle \left\{\xi_\alpha\right\}_{\alpha \in \left\{1,
\ldots , m\right\}} $
is any {\sl local orthonormal normal framefield} on $ M $
in $ \R^{n+m}({\widetilde k})  $.
  One can  remark that
 $ \rho^\perp   =  0 $ if and only if
 $ R^\perp = 0 $, which means that {\sl the normal
 connection is flat}. This follows from
 (\ref{normcurvtens}) and as it was already
  observed by E. Cartan (see \cite{ECA}), is
  equivalent to the simultaneous
  diagonalisability of all shape operators
  $ A_\xi $ of $ M $ in    $ \R^{n+m}({\widetilde k})  $.

\subsection{Wintgen ideal submanifolds satisfying some pseudo-symmetry type curvature conditions}
		\label{ciso2}

	In classical differential geometry, for a surface
$M^2$ in a Euclidean $3$-space $\E^3$, the well-known Euler inequality is
given by
$$
K \leq H^2  ,
$$
whereby $K$ is the intrinsic {\sl Gauss curvature} of $M^2$ and
$H^2$ is the extrinsic
{\sl squared mean curvature} of $M^2$ in  $\E^3$, at once
follows from the fact that $K = k_1k_2$ and $H =\frac{1}{2}(k_1 + k_2)$
whereby $k_1$ and $k_2$ denote
the {\sl principal curvatures}  of $M^2$ in  $\E^3$.
Obviously, $ K = H^2 $ everywhere on $M^2$ if and only the surface
$M^2$ is totally umbilical in $ \E^3$, i.e.
$k_1 = k_2$ at all points of $M^2$, or still, by a
 {\sl theorem of Meunier}, if and only $M^2$
is a part of a {\sl plane} $\E^2$ or of a
{\sl round sphere}  $S^2$ in $\E^3$.
P. Wintgen (see \cite{8}), in the late $19$ seventies, proved that
the {\sl Gauss curvature} $K$ and the
{\sl squared mean curvature} $H^2$ and the
{\sl normal curvature} $K^\perp$ of any
surface $M^2$ in $\E^4$ always
satisfy the inequality
$$
K \leq H^2 - K^\perp
$$
and that actually {\sl the equality} if and only
if the curvature ellipse of $M^2$ in $\E^4$ is a circle.
We recall that the {\sl ellipse of curvature}
at a point $p$ of $M^2$ is defined as
$$ 
\Sigma_p
= \left\{ 
h(X,X) \ \Big\vert \ X \in  T_p M \ \text{and} \ \left\|X\right\|=1 \right\}. 
$$
The {\sl ellipse of curvature} is the analogue of the 	
{\sl Dupin curvature} of an ordinary surface in $\E^3$.

B. Rouxel in \cite{ROUX} and V. Guadalupe and L. Rodriguez in  \cite{GuaRod}  extended
Wintgen's inequality to surfaces of
arbitrary codimension in real space forms
$\R^{2+m}(c)$ with $ m \geq 2$.
Also, B.-Y. Chen extended Wintgen's inequality
in  \cite{{Ch12}, {Ch13}} to surfaces
in pseudo-Euclidean $4$-space $\E^{4}_{2}(c)$ with a neutral metric.

In $1999$, P. J. De Smet, F. Dillen, L. Verstraelen
and L. Vrancken proved in \cite{1} the Wintgen inequality
 $$
   \rho  \leq H^2 - \rho^\perp  + {\widetilde k}
 $$
for all submanifolds $M^n$ of codimension 2 in all real space forms,
whereby $\rho$ is the normalised scalar
curvature of the Riemannian manifold M, and
whereby $H^2$ and $\rho^\perp$, are the squared mean curvature and the normalized normal
scalar curvature of M in the ambient space, respectively, characterizing
the equality in terms of the shape operators of $M^n$ in $\R^{n+2}({\widetilde k})$.
And in \cite{1} they proposed a conjecture of Wintgen inequality for
general Riemannian submanifolds in real space forms,
which was later well-known as the DDVV conjecture.
This conjecture was proven to be true by
Z. Lu (see \cite{2}), and by J. Ge and Z. Tang (see \cite{{4}, {5}})
independently.
In  \cite{BYCH2021}, B.-Y. Chen provided a comprehensive survey on developments 
in Wintgen inequality and Wintgen ideal submanifolds, we also refer to the recent
article of G.-E. Vilcu {\cite[Chapter 7] {Vilcu-2022}}.

The main purpose of the present article is to study
	the so-called {\sl Wintgen ideal submanifolds}.
	The following theorem of J. Ge and Z. Tang
	(see \cite{{4}, {5}})  and Lu (see \cite{2}) states
	 {\sl the Wintgen inequality} for
	submanifolds $ M =  M^n  $ of any  codimension
	$ m \geq 1 $ in a real space form
	$ \R^{n+m}({\widetilde k}) $, and
	{\sl characterizes its equality case}.
   As it is stated in  \cite{4}, it concerns
	{\sl a basic general optimal inequality between}
	likely
	{\sl the most primitive scalar valued geometric quantities}
	that can be defined
	{\sl on submanifolds} as {\sl intrinsic invariant}
	it involves
	{\sl the scalar curvature}
	and as
	{\sl extrinsic invariants}
	it involves
	{\sl the scalar normal curvature}
	and {\sl the squared mean curvature}.

\medskip

  \noindent {\bf Theorem $A$.}
{\rm (see \cite{{3}, {4}, {5}, {2}} )} \label{ChoiLu1}  \  
{\sl   Let $ M = M^n $  be
 a submanifold
  of codimension $ m $
   in a real space form    $ \R^{n+m}({\widetilde k})  $,
	$ n \geq 4$ and $m \geq 2$. Then
   $$
   \mathrm{ (*)}   \quad
	\rho  \leq  H^2 - \rho^\perp  + {\widetilde k} ~;
 $$
  and
  in $ \mathrm{ (*) } $ actually the equality
   holds  if and
  only if,
  with respect to some suitable adapted orthonormal
 frame
  $  \left\{E_i , \xi_\alpha\right\}    $
  on $ M $ in $ \R^{n+m}({\widetilde k})  $, the
  shape operators are given by
\begin{equation}\label{ChoiLu2}
\begin{aligned}
  A_1   &=  \begin{pmatrix}
      a & \mu & 0 & \cdots & 0 \\
    \mu & a & 0 & \cdots & 0 \\
      0 & 0 & a & \cdots & 0 \\
 \vdots & \vdots & \vdots & \ddots & \vdots \\
      0 & 0 & 0 & \cdots & a\\
             \end{pmatrix}  , \  
 A_2     =  \begin{pmatrix}
b + \mu & 0 & 0 & \cdots & 0 \\
      0 & b - \mu & 0 & \cdots & 0 \\
      0 & 0 & b & \cdots & 0 \\
 \vdots & \vdots & \vdots & \ddots & \vdots \\
      0 & 0 & 0 & \cdots & b \\
               \end{pmatrix} , \ 
 A_3     =  \begin{pmatrix}
      c &     0  & 0 & \cdots & 0 \\
   0    & c & 0 & \cdots & 0 \\
   0    & 0 & c & \cdots & 0 \\
 \vdots & \vdots & \vdots & \ddots & \vdots \\
      0 & 0 & 0 & \cdots & c \\
                \end{pmatrix} , \\
A_4 &= \ldots = A_m = 0 .
\end{aligned}
 \end{equation}
 where
 $ a $, $ b $, $ c $ and
 $  \mu $ are real functions on $ M $.
 \ $ \square  $
			}

\medskip

\begin{definition}
The submanifold $ M $ in $ \R^{n+m}({\widetilde k})  $
which satisfy the equality
  $$
   \mathrm{ (\widetilde{*})}   \quad
	\rho  = H^2 - \rho^\perp  + {\widetilde k}
 $$
in the {\sl Wintgen's general inequality $ \mathrm{ (*)} $}
is called  a {\sl Wintgen ideal submanifold}.
In such case,
the frames $  \left\{E_i , \xi_\alpha\right\}    $
in which the shape operators assume the forms
of (\ref{ChoiLu2}) will further be called
the  {\sl Choi-Lu frames},  and
the corresponding tangent
$ E_1E_2 $-planes will be called
the {\sl Choi-Lu planes} of $ M $ in $ \R^{n+m}({\widetilde k})  $.
\end{definition}

\medskip

\noindent  From  \cite{15} (see also \cite{{13}, {12}})
  we recall the following results.

\medskip

  \noindent {\bf Theorem $B$.}  \  
{\sl
 A     Wintgen ideal submanifold
$ M = M^n $ of dimension $ n \geq 4 $ and of codimension $ m $
in a real space form
$  \R^{n+m}({\widetilde k})  $ ($ m \geq 2 $)
 is a Deszcz symmetrical  Riemannian manifold
if and only if   $ M $ is   totally umbilical (with Deszcz sectional curvature
$ L_R = 0 $),   or,   $ M $ is
a minimal or  pseudo-umbilical
submanifold (with  $ L_R =  {\widetilde k} + H^2 $) of this space form
$   \R^{n+m}({\widetilde k})  $.
	\     $ \square $

     }

\medskip

\noindent {\bf Theorem $C$.}  \  
{\sl  A     Wintgen ideal submanifold
$ M = M^n $ of dimension $ n \geq 4 $ and of codimension $m$ in a real space form
  $  \R^{n+m}({\widetilde k})  $
 ($ m \geq 2 $)  is Deszcz symmetric
if and only if
$M$  is  Deszcz Ricci-symmetric.
	\     $ \square $
      }

\medskip

\noindent {\bf Theorem $D$.}  \  
{\sl   
	Let $ M = M^n $ be
 a   Wintgen ideal submanifold
 of dimension $ n \geq 4 $ and of codimension $ m \geq 2 $ in a real space form
$  \R^{n+m}({\widetilde k}) $.
  Then
$	 M  $ is a Riemannian manifold with
 a pseudo-symmetric  conformal Weyl tensor  $C$.
		 }

\medskip

\noindent {\bf Theorem $E$.}  \
{\sl   
	Let $ M = M^n $ be
 a   Wintgen ideal submanifold
 of dimension $ n \geq 4 $ and of codimension $ m \geq 2 $ in a real space form
$  \R^{n+m}({\widetilde k})  $.

\noindent 	(i)   $M$ is conformally flat if and only if $M$ is
	a totally umbilical submanifold in  $  \R^{n+m}({\widetilde k}) $
	(and, hence, $M$ is a real space form).

	\noindent (ii)	If $M$ is not a conformally flat submanifold,
	then 	$ M  $ has  a pseudo-symmetric conformal Weyl tensor $C$
and the corresponding function   is
	given by
	$$
	L_C  =   \dfrac{n-3}{(n-1)(n-2)}\Big(\kappa - n(n-1)\inf K\Big) .
	$$	
  where $\kappa$ is the scalar curvature of $M$.    $     \square		 $}

We also refer to \cite{{14}, {MPTAP}, {SEB}, {SEBMPTAP}, {SEN}}
for further results on Wintgen ideal submanifolds.

\section{Main results on Wintgen ideal submanifolds}

  \subsection{Main results}

We consider  mainly
the
Wintgen ideal
submanifold
$ M = M^n $ of dimension $ n > 3 $
and of codimension $ m \geq 2 $ in a real space form
$  \R^{n+m}({\widetilde k})  $
 and the Theorem $A$ in the last section  \ref{ciso2}.
 We  investigate $M$
 satisfying one of the following pseudo-symmetry type curvature conditions:
\begin{enumerate}
	\item[(i)] the tensor $ R \cdot C  $ and the Tachibana tensor  $ Q(g , R) $
  (resp., $ Q(g , C) $,  $ Q(g , g\wedge\mbox{\rm Ricc}) $,
	   $ Q(\mbox{\rm Ricc} , R) $, $Q(\mbox{\rm Ricc} , C) $)
	   are linearly dependent;

	\item[(ii)] the tensor $   C \cdot R $ and the Tachibana tensor  $ Q(g , R) $
	(resp., $ Q(g , C) $,  $ Q(g , g\wedge\mbox{\rm Ricc}) $,
	   $ Q(\mbox{\rm Ricc} , R) $, $Q(\mbox{\rm Ricc} , C) $)
	   are linearly dependent;
				
	\item[(iii)] the tensor $ R \cdot C - C \cdot R $ nd the Tachibana tensor
	$ Q(g , R) $
	(resp., $ Q(g , C) $,  $ Q(g , g\wedge\mbox{\rm Ricc}) $,
	   $ Q(\mbox{\rm Ricc} , R) $, $Q(\mbox{\rm Ricc} , C) $)
	   are linearly dependent. 	 	
\end{enumerate}

\begin{theoreme}
  Let $ M= M^n $  be
 a  Wintgen ideal
  submanifold
  of codimension $ m $
   in a real space form
	$ \R^{n+m}({\widetilde k}) $,
	$ n \geq 4 $ and  $ m \geq 2 $.
	If $ {\widetilde k} > 0 $,
	then 	 $M$ is Weyl-semi-symmetric
   	if and only if
  $M$ is totally umbilical in $ \R^{n+m}({\widetilde k}) $
  (and hence $M$ is conformally flat). \     $ \square$
 \end{theoreme}

\begin{theoreme}
  Let $ M = M^n $  be  a  Wintgen ideal
  submanifold of codimension $ m $
   in a real space form
		$  \R^{n+m}({\widetilde k}) $,
	where    $ n \geq 4 $ and  $ m \geq 2 $.
	If  ${\widetilde k} \leq 0 $,
	then		 $M$ is Weyl-semi-symmetric
   	if and only if
	 (i) either $ M$ is totally umbilical
 in  $   \R^{n+m}({\widetilde k})   $
(and hence  $M$ is conformally flat),
(ii) or
 with respect to some suitable adapted
Choi-Lu frames
  $  \left\{E_i , \xi_\alpha\right\}    $
  on $ M $ in $ \R^{n+m}({\widetilde k})  $, the
  shape operators are given by
\begin{equation}\label{ChoiLu611935}
\begin{aligned}
  A_1 &=  \begin{pmatrix}
  0   & \mu & 0 & \cdots & 0 \\
  \mu &  0  & 0 & \cdots & 0 \\
  0   &  0  & 0 & \cdots & 0 \\
  \vdots & \vdots & \vdots & \ddots & \vdots \\
  0 & 0 & 0 & \cdots & 0\\
  \end{pmatrix}  , \ 
 A_2  =  \begin{pmatrix}
    \mu & 0 & 0 & \cdots & 0 \\
  0 &   -\mu & 0 & \cdots & 0 \\
  0 & 0 & 0 & \cdots & 0 \\
  \vdots & \vdots & \vdots & \ddots & \vdots \\
  0 & 0 & 0 & \cdots & 0 \\
  \end{pmatrix}
,  \   
 A_3  =  \begin{pmatrix}
  c &     0  & 0 & \cdots & 0 \\
  0      & c & 0 & \cdots & 0 \\
  0 & 0 & c & \cdots & 0 \\
  \vdots & \vdots & \vdots & \ddots & \vdots \\
  0 & 0 & 0 & \cdots & c \\
  \end{pmatrix}
, \     \\
A_4  &= \ldots = A_m = 0 ,
\end{aligned}
 \end{equation}
 where
$c$,   $  \mu $ are real functions on $ M $ such that $ \mu \ne 0 $
and $ c^2 = -{\widetilde k} $.
	In this second case, $M$ is
	a minimal
 or    pseudo-umbilical submanifold
	in $ \R^{n+m}({\widetilde k}) $  . \ $\square$

  \end{theoreme}

\begin{corollaire}
  Let $ M = M^n $  be  a  Wintgen ideal
  submanifold of codimension $ m $
   in a real space form
		$  \R^{n+m}({\widetilde k})  $,
	where    $ n \geq 4 $ and  $ m \geq 2 $.
	If  ${\widetilde k} = 0 $,
	then	 	 $M$ is Weyl-semi-symmetric
    	if and only if
(i) either $ M$ is totally umbilical
 in  $   \R^{n+m}    $
(and hence  $M$ is conformally flat),
(ii) or   with respect to some suitable adapted
Choi-Lu frames
  $  \left\{E_i , \xi_\alpha\right\}    $
  on $ M $ in $ \R^{n+m}   $, the
  shape operators are given by
\begin{equation}\label{ChoiLu61139999}
\begin{aligned}
  A_1 &=  \begin{pmatrix}
  0   & \mu & 0 & \cdots & 0 \\
  \mu &  0  & 0 & \cdots & 0 \\
  0   &  0  & 0 & \cdots & 0 \\
  \vdots & \vdots & \vdots & \ddots & \vdots \\
  0 & 0 & 0 & \cdots & 0\\
  \end{pmatrix}  , \   
 A_2  =  \begin{pmatrix}
    \mu & 0 & 0 & \cdots & 0 \\
  0 &   -\mu & 0 & \cdots & 0 \\
  0 & 0 & 0 & \cdots & 0 \\
  \vdots & \vdots & \vdots & \ddots & \vdots \\
  0 & 0 & 0 & \cdots & 0 \\
  \end{pmatrix}
,   \\
 A_3 &= A_4 = \ldots = A_m = 0 ,
\end{aligned}
 \end{equation}
 where $\mu $ is a real function  on $ M $
 such that $ \mu \ne 0 $.
 In this second case,
  $M$ is minimal in $ \R^{n+m}   $.  \ $\square$

   \end{corollaire}


\begin{theoreme}
  Let $ M= M^n $  be
 a  Wintgen ideal
  submanifold
  of codimension $ m $
   in a real space form
	$ \R^{n+m}({\widetilde k}) $,
	$ n \geq 4 $ and  $ m \geq 2 $.
	Then 	
  	$$ C \cdot R = 0  $$
  	if and only if
  $M$ is totally umbilical in $ \R^{n+m}({\widetilde k}) $
  (and hence $M$ is conformally flat). \     $ \square$
 \end{theoreme}


\begin{theoreme}
  Let $ M= M^n $  be
 a  Wintgen ideal
  submanifold
  of codimension $ m $
   in a real space form
	$ \R^{n+m}({\widetilde k}) $,
	$ n \geq 4 $ and  $ m \geq 2 $.
	Then 	
  	$$  R \cdot C  -   C \cdot R = 0  $$
  	if and only if
  $M$ is totally umbilical in $ \R^{n+m}({\widetilde k}) $
  (and hence $M$ is conformally flat). \     $ \square$
 \end{theoreme}


\begin{theoreme}
  Let $ M= M^n $  be
 a  Wintgen ideal
  submanifold
  of codimension $ m $
   in a real space form
	$ \R^{n+m}({\widetilde k}) $,
	$ n \geq 4 $ and  $ m \geq 2 $.
	If $ {\widetilde k} > 0 $,
	then 	
  the tensors  	$ R \cdot C   $ and
		$  Q(g , R) $   are   linearly dependent
if and only if $M$ is totally umbilical in
$ \R^{n+m}({\widetilde k}) $
(and hence 	$M$ is conformally flat). \     $ \square$
 \end{theoreme}


\begin{theoreme}
  Let $ M = M^n $  be  a  Wintgen ideal
  submanifold of codimension $ m $
   in a real space form
		$  \R^{n+m}({\widetilde k}) $,
	where    $ n \geq 4 $ and  $ m \geq 2 $.
	If  ${\widetilde k} \leq 0 $,
	then	  the tensors  $ R \cdot C   $
    and
		$  Q(g , R) $
	 are linearly dependent	if and only if
	 (i) either $ M$ is totally umbilical
 in  $   \R^{n+m}({\widetilde k})   $
(and hence  $M$ is conformally flat),
(ii) or
 with respect to some suitable adapted
Choi-Lu frames
  $  \left\{E_i , \xi_\alpha\right\}    $
  on $ M $ in $ \R^{n+m}({\widetilde k})  $, the
  shape operators are given by
\begin{equation}\label{ChoiLu611}
\begin{aligned}
  A_1 &=  \begin{pmatrix}
  0   & \mu & 0 & \cdots & 0 \\
  \mu &  0  & 0 & \cdots & 0 \\
  0   &  0  & 0 & \cdots & 0 \\
  \vdots & \vdots & \vdots & \ddots & \vdots \\
  0 & 0 & 0 & \cdots & 0\\
  \end{pmatrix}  , \ 
 A_2  =  \begin{pmatrix}
    \mu & 0 & 0 & \cdots & 0 \\
  0 &   -\mu & 0 & \cdots & 0 \\
  0 & 0 & 0 & \cdots & 0 \\
  \vdots & \vdots & \vdots & \ddots & \vdots \\
  0 & 0 & 0 & \cdots & 0 \\
  \end{pmatrix}
, \  
 A_3  =  \begin{pmatrix}
  c &     0  & 0 & \cdots & 0 \\
  0      & c & 0 & \cdots & 0 \\
  0 & 0 & c & \cdots & 0 \\
  \vdots & \vdots & \vdots & \ddots & \vdots \\
  0 & 0 & 0 & \cdots & c \\
  \end{pmatrix}
,   \\
A_4  &= \ldots = A_m = 0 ,
\end{aligned}
 \end{equation}
 where
$c$,   $  \mu $ are real functions on $ M $ such that $ \mu \ne 0 $
and $ c^2 = -{\widetilde k} $.
	Moreover, $M$ is
	Weyl-semi-symmetric~;
	and,   in the second case,
	$M$ is
	a minimal  or  pseudo-umbilical submanifold
	in $ \R^{n+m}({\widetilde k}) $.  \  $\square$
	
  \end{theoreme}

\begin{corollaire}
  Let $ M = M^n $  be  a  Wintgen ideal
  submanifold of codimension $ m $
   in a real space form
		$  \R^{n+m}({\widetilde k})  $,
	where    $ n \geq 4 $ and  $ m \geq 2 $.
	If  ${\widetilde k} = 0 $,
	then	  the tensors  $ R \cdot C   $
    and
		$  Q(g ,   R) $
	 are linearly dependent	if and only if
	 (i) either $ M$ is totally umbilical
 in  $   \R^{n+m}    $
(and hence  $M$ is conformally flat),
(ii) or   with respect to some suitable adapted
Choi-Lu frames
  $  \left\{E_i , \xi_\alpha\right\}    $
  on $ M $ in $ \R^{n+m}   $, the
  shape operators are given by
\begin{equation}\label{ChoiLu6113}
\begin{aligned}
  A_1 &=  \begin{pmatrix}
  0   & \mu & 0 & \cdots & 0 \\
  \mu &  0  & 0 & \cdots & 0 \\
  0   &  0  & 0 & \cdots & 0 \\
  \vdots & \vdots & \vdots & \ddots & \vdots \\
  0 & 0 & 0 & \cdots & 0\\
  \end{pmatrix}  , \   
 A_2  =  \begin{pmatrix}
    \mu & 0 & 0 & \cdots & 0 \\
  0 &   -\mu & 0 & \cdots & 0 \\
  0 & 0 & 0 & \cdots & 0 \\
  \vdots & \vdots & \vdots & \ddots & \vdots \\
  0 & 0 & 0 & \cdots & 0 \\
  \end{pmatrix}
, \\
 A_3 &= A_4 = \ldots = A_m = 0 ,
\end{aligned}
 \end{equation}
 where $\mu $ is a real function  on $ M $
 such that $ \mu \ne 0 $.
 In this second case,
  $M$ is minimal in $ \R^{n+m}   $.  \ $\square$

   \end{corollaire}


   \begin{theoreme}  \label{cr.LQgCorigin61}   
  Let $ M = M^n$ be a Wintgen ideal
	submanifold of codimension
  $ m $ in a real space form
		$  \R^{n+m}({\widetilde k}) $,
  $ n \geq 4 $ and  $ m \geq 1 $. Then   the tensors
	$ C \cdot R $ and $Q(g, R)$
	 are linearly dependent
if and only if $M$ is totally   umbilical in
 $  \R^{n+m}({\widetilde k}) $ (and hence $M$	
 is conformally flat).  \     $ \square$
       \end{theoreme}


\begin{theoreme}
  Let $ M = M^n $  be
 a  Wintgen ideal
  submanifold
  of codimension $ m $
   in a real space form    $ \R^{n+m}({\widetilde k})  $,
	$ n \geq 4 $ and  $ m \geq 2 $.
	If $ {\widetilde k} > 0 $,
	then
  the tensors  $ R \cdot C - C \cdot R $
    and
		$  Q(g ,   R) $     are   linearly dependent
if and only if $M$ is totally umbilical in
$ \R^{n+m}({\widetilde k}) $
(and hence 	$M$ is conformally flat). \      $ \square$
 \end{theoreme}

\begin{theoreme}
  Let $ M = M^n$  be
 a  Wintgen ideal
  submanifold
  of codimension $ m $
   in a real space form
	$ \R^{n+m}({\widetilde k})  $,
	where       $ n \geq 4 $ and  $ m \geq 2 $.
	If ${\widetilde k} \leq 0 $,
	then	  the tensors
  $ R \cdot C - C \cdot R $
    and
		$  Q(g ,   R) $  are linearly dependent
		if and only if
 (i) either $ M$ is totally umbilical
 in  $   \R^{n+m}({\widetilde k})   $
(and hence  $M$ is conformally flat),
(ii) or   with respect to some suitable adapted
Choi-Lu frames
  $  \left\{E_i , \xi_\alpha\right\}    $
  on $ M $ in $ \R^{n+m}({\widetilde k})  $, the
  shape operators are given by
\begin{equation}\label{ChoiLu6112}
\begin{aligned}
  A_1 &=  \begin{pmatrix}
  0   & \mu & 0 & \cdots & 0 \\
  \mu &  0  & 0 & \cdots & 0 \\
  0   &  0  & 0 & \cdots & 0 \\
  \vdots & \vdots & \vdots & \ddots & \vdots \\
  0 & 0 & 0 & \cdots & 0\\
  \end{pmatrix}  , \   
 A_2  =  \begin{pmatrix}
    \mu & 0 & 0 & \cdots & 0 \\
  0 &   -\mu & 0 & \cdots & 0 \\
  0 & 0 & 0 & \cdots & 0 \\
  \vdots & \vdots & \vdots & \ddots & \vdots \\
  0 & 0 & 0 & \cdots & 0 \\
  \end{pmatrix}
, \  
 A_3 =  \begin{pmatrix}
  c &     0  & 0 & \cdots & 0 \\
  0      & c & 0 & \cdots & 0 \\
  0 & 0 & c & \cdots & 0 \\
  \vdots & \vdots & \vdots & \ddots & \vdots \\
  0 & 0 & 0 & \cdots & c \\
  \end{pmatrix}
,    \\
A_4  &= \ldots = A_m = 0 ,
\end{aligned}
 \end{equation}
 where
$c$,   $  \mu $ are real functions on $ M $ such that $ \mu \ne 0 $
and $ c^2 = -{\widetilde k} $.
  Moreover,
	    $$
     R \cdot C - C \cdot R
     = -\frac{2(n-3)\mu^2}{(n-1)(n-2)}Q(g, R) ~; 
   $$
	and, in the second case,
		$M$ is
   a minimal  or pseudo-umbilical submanifold
	in $ \R^{n+m}({\widetilde k}) $.  \      $\square  $
   \end{theoreme}

\begin{corollaire}
  Let $ M = M^n$  be
 a  Wintgen ideal
  submanifold
  of codimension $ m $
   in a real space form
	$ \R^{n+m}({\widetilde k})  $,
	where       $ n \geq 4 $ and  $ m \geq 2 $.
	If ${\widetilde k} = 0 $,
	then	  the tensors
  $ R \cdot C - C \cdot R $
    and
		$  Q(g , R) $  are linearly dependent
		if and only if
 (i) either $ M$ is totally umbilical
 in  $   \R^{n+m}  $
(and hence  $M$ is conformally flat),
(ii) or  with respect to some suitable adapted
Choi-Lu frames
  $  \left\{E_i , \xi_\alpha\right\}    $
  on $ M $ in $ \R^{n+m}   $, the
  shape operators are given by
\begin{equation}\label{ChoiLu61123}
\begin{aligned}
  A_1 &=  \begin{pmatrix}
  0   & \mu & 0 & \cdots & 0 \\
  \mu &  0  & 0 & \cdots & 0 \\
  0   &  0  & 0 & \cdots & 0 \\
  \vdots & \vdots & \vdots & \ddots & \vdots \\
  0 & 0 & 0 & \cdots & 0\\
  \end{pmatrix}  , \  
 A_2  =  \begin{pmatrix}
    \mu & 0 & 0 & \cdots & 0 \\
  0 &   -\mu & 0 & \cdots & 0 \\
  0 & 0 & 0 & \cdots & 0 \\
  \vdots & \vdots & \vdots & \ddots & \vdots \\
  0 & 0 & 0 & \cdots & 0 \\
  \end{pmatrix}
,    \\
 A_3 &= A_4 = \ldots = A_m = 0 ,
\end{aligned}
 \end{equation}
 where
  $  \mu $ is a  real function  on $ M $
 such that $ \mu \ne 0 $.
 In this second case,  $M$ is minimal
 in $ \R^{n+m}   $. \      $\square  $
   \end{corollaire}


   \begin{theoreme}  \label{cr.LQgCorigin3}   
  Let $ M = M^n $ be  a Wintgen ideal
	submanifold of codimension
  $ m $ in a real space form
	$ \R^{n+m}({\widetilde k})  $,
  $ n \geq 4 $ and  $ m \geq 1 $. Then
	 the tensors   $ R \cdot C $ 	 and
              $	Q(g, C)$
  are linearly dependent  if and only if
 (i) either $ M$  is totally umbilical
 in $   \R^{n+m}({\widetilde k})   $
(and hence  $M$  is conformally flat),
(ii) or   with respect to some suitable adapted
Choi-Lu frames
   $  \left\{E_i , \xi_\alpha\right\}    $
  on $ M $ in $ \R^{n+m}({\widetilde k})   $, the
  shape operators are given by
\begin{equation}\label{ChoiLu23}
\begin{aligned}
  A_1 &=  \begin{pmatrix}
  0   & \mu & 0 & \cdots & 0 \\
  \mu &  0  & 0 & \cdots & 0 \\
  0   &  0  & 0 & \cdots & 0 \\
  \vdots & \vdots & \vdots & \ddots & \vdots \\
  0 & 0 & 0 & \cdots & 0\\
  \end{pmatrix}  , \  
 A_2  =  \begin{pmatrix}
    \mu & 0 & 0 & \cdots & 0 \\
  0 &   -\mu & 0 & \cdots & 0 \\
  0 & 0 & 0 & \cdots & 0 \\
  \vdots & \vdots & \vdots & \ddots & \vdots \\
  0 & 0 & 0 & \cdots & 0 \\
  \end{pmatrix}
, \   
 A_3  =  \begin{pmatrix}
  c &     0  & 0 & \cdots & 0 \\
  0      & c & 0 & \cdots & 0 \\
  0 & 0 & c & \cdots & 0 \\
  \vdots & \vdots & \vdots & \ddots & \vdots \\
  0 & 0 & 0 & \cdots & c \\
  \end{pmatrix}
,   \\
A_4  &= \ldots = A_m = 0 ,
\end{aligned}
 \end{equation}
 where
$ c $ and
 $  \mu $ are real functions on $ M $ such that $ \mu \ne 0 $.
Moreover,
$$
 R \cdot C =  ({\widetilde k} + c^2) Q(g, C)  ~;
 $$
and, in the second case,
  $M$ is a minimal or pseudo-umbilical
submanifold  in $ \R^{n+m}(\widetilde{k}) $  . \      $\square$
\end{theoreme}

   \begin{theoreme}  \label{cr.LQgCorigin1}
  Let $ M = M^n$ be a Wintgen ideal
	submanifold of codimension
  $ m $ in a real space form
	$ \R^{n+m}({\widetilde k})  $,
  $ n \geq 4 $ and  $ m \geq 1 $. Then   the tensors
	$  C \cdot R  $
	  and
		$ Q(g, C)$  are linearly dependent
 if and only if $ M$ is totally umbilical
in  $   \R^{n+m}({\widetilde k})   $.
In this case, $M$ is conformally flat. \    $ \square$
     \end{theoreme}

  \begin{theoreme}    \label{r.c.LQgRLQSR51}\  
  Let $ M = M^n$  be  a Wintgen ideal
	submanifold of codimension
  $ m $ in a real space form
	$\R^{n+m}({\widetilde k})  $,
  $ n \geq 4 $ and  $ m \geq 1 $. Then   the tensors
 $R\cdot C -   C \cdot R$
	  and
$Q(g, C)$  are linearly dependent
 if and if
$ M $  is totally umbilical
in  $   \R^{n+m}({\widetilde k})   $.
 In this  case, $M$ is  conformally flat.
\     $ \square $
     \end{theoreme}


  \begin{theoreme}   \label{cr.LQgCorigin6162}
  Let $ M = M^n $  be  a Wintgen ideal
	submanifold of codimension
  $ m $ in a real space form
	$ \R^{n+m}({\widetilde k})  $,
  $ n \geq 4 $ and  $ m \geq 1 $.
	If $ {\widetilde k} > 0$,
	then
  the tensors  $ R \cdot C $ and $Q(g, g\wedge\mbox{\rm Ricc})$
	 are   linearly dependent
if and only if $M$ is totally umbilical in
$ \R^{n+m}({\widetilde k}) $
(and hence 	$M$ is conformally flat). \     $ \square$
			\end{theoreme}

   \begin{theoreme}  
	           \label{cr.LQgCorigin6163}   
  Let $ M = M^n$  be  a Wintgen ideal
	submanifold of codimension
  $ m $ in a real space form
	$  \R^{n+m}({\widetilde k})  $,
  $ n \geq 4 $ and  $ m \geq 1 $.
   If $ {\widetilde k} \leq 0$,
	then   the tensors
	$ R \cdot C $ and $Q(g, g\wedge\mbox{\rm Ricc})$
		 are linearly dependent	if and only if
			(i) either $ M$  is totally umbilical
in  $   \R^{n+m}({\widetilde k})   $
(and hence  $M$ is conformally flat),
(ii) or   with respect to some suitable adapted
Choi-Lu frames
  $  \left\{E_i , \xi_\alpha\right\}    $
  on $ M $ in $ \R^{n+m}({\widetilde k})  $, the
  shape operators are given by
\begin{equation}\label{ChoiLu6163}
\begin{aligned}
  A_1 &=  \begin{pmatrix}
  0   & \mu & 0 & \cdots & 0 \\
  \mu &  0  & 0 & \cdots & 0 \\
  0   &  0  & 0 & \cdots & 0 \\
  \vdots & \vdots & \vdots & \ddots & \vdots \\
  0 & 0 & 0 & \cdots & 0\\
  \end{pmatrix}  ,   \ 
 A_2  =  \begin{pmatrix}
    \mu & 0 & 0 & \cdots & 0 \\
  0 &   -\mu & 0 & \cdots & 0 \\
  0 & 0 & 0 & \cdots & 0 \\
  \vdots & \vdots & \vdots & \ddots & \vdots \\
  0 & 0 & 0 & \cdots & 0 \\
  \end{pmatrix}
, \    
 A_3  =  \begin{pmatrix}
  c     & 0      & 0      & \cdots & 0 \\
  0     & c      & 0      & \cdots & 0 \\
  0     & 0      & c      & \cdots & 0 \\
 \vdots & \vdots & \vdots & \ddots & \vdots \\
  0     & 0      & 0      & \cdots & c \\
  \end{pmatrix}
,   \\
A_4  &= \ldots = A_m = 0 ,
\end{aligned}
 \end{equation}
 where
$ c $ and
 $  \mu $ are real functions on $ M $ such that $ \mu \ne 0 $
and $ c^2 = - {\widetilde k} $.
Moreover, $ M$ is Weyl-semi-symmetric~;
and,
in the second case, $M$ is
	 a minimal or
  pseudo-umbilical
	submanifold  in $ \R^{n+m}({\widetilde k}) $. \  $\square$

      \end{theoreme}


   \begin{corollaire}  
	           \label{cr.LQgCorigin616310}
  Let $ M = M^n$  be  a Wintgen ideal
	submanifold of codimension
  $ m $ in a real space form
	$  \R^{n+m}({\widetilde k})  $,
  $ n \geq 4 $ and  $ m \geq 1 $.
   If $ {\widetilde k} = 0$,
	then   the tensors
	$ R \cdot C $ and $Q(g, g\wedge\mbox{\rm Ricc})$
		 are linearly dependent	if and only if
			(i) either $ M$  is totally umbilical
in  $   \R^{n+m}  $
(and hence  $M$ is conformally flat),
(ii) or   with respect to some suitable adapted
Choi-Lu frames
  $  \left\{E_i , \xi_\alpha\right\}    $
  on $ M $ in $ \R^{n+m}  $, the
  shape operators are given by
\begin{equation}\label{ChoiLu616359}
\begin{aligned}
  A_1 &=  \begin{pmatrix}
  0   & \mu & 0 & \cdots & 0 \\
  \mu &  0  & 0 & \cdots & 0 \\
  0   &  0  & 0 & \cdots & 0 \\
  \vdots & \vdots & \vdots & \ddots & \vdots \\
  0 & 0 & 0 & \cdots & 0\\
  \end{pmatrix}  , \   
 A_2  =  \begin{pmatrix}
    \mu & 0 & 0 & \cdots & 0 \\
  0 &   -\mu & 0 & \cdots & 0 \\
  0 & 0 & 0 & \cdots & 0 \\
  \vdots & \vdots & \vdots & \ddots & \vdots \\
  0 & 0 & 0 & \cdots & 0 \\
  \end{pmatrix}
, \\
A_3 &= A_4 \ldots = A_m = 0 ,
\end{aligned}
 \end{equation}
 where
  $  \mu $ is a real function  on $ M $
  such that $ \mu \ne 0 $.
 In this second case,
 $M$ is   minimal in $ \R^{n+m}  $.
   \  $\square$

      \end{corollaire}


   \begin{theoreme}   \label{cr.LQgCorigin6161}   
  Let $ M = M^n $ be  a Wintgen ideal
	submanifold of codimension
  $ m $ in a real space form
	$ \R^{n+m}({\widetilde k})  $,
  $ n \geq 4 $ and  $ m \geq 1 $. Then
	$ C \cdot R $ and  $Q(g, g\wedge\mbox{\rm Ricc})$
	are linearly dependent		if and only if   $ M$  is totally umbilical
in  $   \R^{n+m}({\widetilde k})   $
(and hence  $M$ is conformally flat). \  $\square$
       \end{theoreme}

  \begin{theoreme} 
	           \label{cr.LQgGwGorigin6171}
  Let $ M = M^n $  be  a Wintgen ideal
	submanifold of codimension
  $ m $ in a real space form
	$ \R^{n+m}({\widetilde k})  $,
  $ n \geq 4 $ and  $ m \geq 1 $. Then   the tensors
	$ R \cdot C - C \cdot R $ and $Q(g, g\wedge\mbox{\rm Ricc})$
		 are linearly dependent
 $ M$  is totally umbilical
in  $   \R^{n+m}({\widetilde k})   $
(and hence  $M$ is conformally flat). \  $\square$
\end{theoreme}


  \begin{theoreme}   \label{cr.LQgCorigin6183}   
  Let $ M = M^n $  be  a Wintgen ideal
	submanifold of codimension
  $ m $ in a real space form
	$ \R^{n+m}({\widetilde k})  $,
  $ n \geq 4 $ and  $ m \geq 1 $.
	If $ {\widetilde k} > 0$,
then   the tensors  $ R \cdot C $ and   $Q(\mbox{\rm Ricc} , R)$
	 	 are   linearly dependent
if and only if $M$ is totally umbilical in
$ \R^{n+m}({\widetilde k}) $
(and hence 	$M$ is conformally flat). \     $ \square$
       \end{theoreme}

   \begin{theoreme}   \label{cr.LQgCorigin618491}   
  Let $ M = M^n $  be  a Wintgen ideal
	submanifold of codimension
  $ m $ in a real space form
	$  \R^{n+m}({\widetilde k})  $,
  $ n \geq 4 $ and  $ m \geq 1 $.
	If $ {\widetilde k} \leq 0$,
	then   the tensors	$ R \cdot C $ and Q(\mbox{\rm Ricc},R)
	are linearly dependent
			if and only if (i) either $ M$  is totally umbilical
in  $   \R^{n+m}({\widetilde k})   $
(and hence  $M$ is conformally flat),
(ii) or   with respect to some suitable adapted
Choi-Lu frames
  $  \left\{E_i , \xi_\alpha\right\}    $
  on $ M $ in $ \R^{n+m}({\widetilde k})  $, the
  shape operators are given by
\begin{equation}\label{ChoiLu6184}
\begin{aligned}
  A_1 &=  \begin{pmatrix}
  0   & \mu & 0 & \cdots & 0 \\
  \mu &  0  & 0 & \cdots & 0 \\
  0   &  0  & 0 & \cdots & 0 \\
  \vdots & \vdots & \vdots & \ddots & \vdots \\
  0 & 0 & 0 & \cdots & 0\\
  \end{pmatrix}  , \    
 A_2  =  \begin{pmatrix}
    \mu & 0 & 0 & \cdots & 0 \\
  0 &   -\mu & 0 & \cdots & 0 \\
  0 & 0 & 0 & \cdots & 0 \\
  \vdots & \vdots & \vdots & \ddots & \vdots \\
  0 & 0 & 0 & \cdots & 0 \\
  \end{pmatrix}
, \   
 A_3  =  \begin{pmatrix}
  c      &   0    & 0      & \cdots & 0 \\
  0      & c      & 0      & \cdots & 0 \\
  0      & 0      & c      & \cdots & 0 \\
  \vdots & \vdots & \vdots & \ddots & \vdots \\
  0      & 0      & 0      & \cdots & c \\
  \end{pmatrix}
,    \\
A_4  &= \ldots = A_m = 0 ,
\end{aligned}
 \end{equation}
 where
$ c $ and
 $  \mu $ are real functions on $ M $ such that $ \mu \ne 0 $
and $ c^2 = - {\widetilde k} $.
Moreover, $M$ is Weyl-semi-symmetric~;
and,
in the second   case,
$M$ is
a minimal or
  pseudo-umbilical submanifold
   in $ \R^{n+m}({\widetilde k})$. \   $\square$
      \end{theoreme}

   \begin{corollaire}   \label{cr.LQgCorigin6188418}   
  Let $ M $  be  a Wintgen ideal
	submanifold of codimension
  $ m $ in a real space form    $   \R^{n+m}({\widetilde k})  $,
  $ n \geq 4 $ and  $ m \geq 1 $.
	  If $ {\widetilde k} = 0$,
	then  the tensors
	$ R \cdot C $ and Q(\mbox{\rm Ricc},R)
	 are linearly dependent
		if and only if (i) either $ M$  is totally umbilical
in  $   \R^{n+m}  $
(and hence  $M$ is conformally flat),
(ii) or   with respect to some suitable adapted
Choi-Lu frames
  $  \left\{E_i , \xi_\alpha\right\}    $
  on $ M $ in $ \R^{n+m} $, the
  shape operators are given by
\begin{equation}\label{ChoiLu618841}
\begin{aligned}
  A_1 &=  \begin{pmatrix}
  0   & \mu & 0 & \cdots & 0 \\
  \mu &  0  & 0 & \cdots & 0 \\
  0   &  0  & 0 & \cdots & 0 \\
  \vdots & \vdots & \vdots & \ddots & \vdots \\
  0 & 0 & 0 & \cdots & 0\\
  \end{pmatrix}  , \   
 A_2  =  \begin{pmatrix}
    \mu & 0 & 0 & \cdots & 0 \\
  0 &   -\mu & 0 & \cdots & 0 \\
  0 & 0 & 0 & \cdots & 0 \\
  \vdots & \vdots & \vdots & \ddots & \vdots \\
  0 & 0 & 0 & \cdots & 0 \\
  \end{pmatrix}
, \\
 A_3 &=  A_4  = \ldots = A_m = 0 ,
\end{aligned}
 \end{equation}
 where
 $  \mu $ is a real function  on $ M $
  such that $ \mu \ne 0 $.
 In this second case,  $M$ is minimal
  in $ \R^{n+m} $. \ $\square$
        \end{corollaire}

  \begin{theoreme}   \label{cr.LQgCorigin6182}   
  Let $ M = M^n$  be  a Wintgen ideal
	submanifold of codimension
  $ m $ in a real space form
	$ \R^{n+m}({\widetilde k})  $,
  $ n \geq 4 $ and  $ m \geq 1 $.
Then   the tensors
	$ C \cdot R $ and $Q(\mbox{\rm Ricc}, R)$
 are linearly dependent	if and only if
	$ M$  is totally umbilical
in  $   \R^{n+m}({\widetilde k})   $
(and hence  $M$ is conformally flat).
 \    $  \square $
       \end{theoreme}

  \begin{theoreme}   \label{cr.LQgCorigin61981}   
  Let $ M = M^n $  be  a Wintgen ideal
	submanifold of codimension
  $ m $ in a real space form
	$ \R^{n+m}({\widetilde k})  $,
  $ n \geq 4 $ and  $ m \geq 1 $.
	Then   the tensors
	$ R \cdot C  - C \cdot R$ and $Q(\mbox{\rm Ricc} , R)$
	 are linearly dependent
if and only if
$ M$  is totally umbilical
in  $   \R^{n+m}({\widetilde k})  $
(and hence  $M$ is conformally flat). \ $\square$
        \end{theoreme}




  \begin{theoreme}   \label{cr.LQgCorigin61883}   
  Let $ M = M^n $  be  a Wintgen ideal
	submanifold of codimension
  $ m $ in a real space form
	$ \R^{n+m}({\widetilde k})  $,
  $ n \geq 4 $ and  $ m \geq 1 $.
	If $ {\widetilde k} > 0$,
	then
	 the tensors  $ R \cdot C $ and   $Q(\mbox{\rm Ricc},C)$
	 are   linearly dependent
if and only if $M$ is totally umbilical in
$ \R^{n+m}({\widetilde k}) $
(and hence 	$M$ is conformally flat). \     $ \square$
       \end{theoreme}

   \begin{theoreme}   \label{cr.LQgCorigin61884}   
  Let $ M = M^n $  be  a Wintgen ideal
	submanifold of codimension
  $ m $ in a real space form
	$  \R^{n+m}({\widetilde k})  $,
  $ n \geq 4 $ and  $ m \geq 1 $.
	If $ {\widetilde k} \leq 0$,
	then   the tensors
	$ R \cdot C $ and Q(\mbox{\rm Ricc},C)
		 are linearly dependent
			if and only if (i) either $ M$  is totally umbilical
in  $   \R^{n+m}({\widetilde k})   $
(and hence  $M$ is conformally flat),
(ii) or   with respect to some suitable adapted
Choi-Lu frames
  $  \left\{E_i , \xi_\alpha\right\}    $
  on $ M $ in $ \R^{n+m}({\widetilde k})  $, the
  shape operators are given by
\begin{equation}\label{ChoiLu618412}
\begin{aligned}
  A_1 &=  \begin{pmatrix}
  0   & \mu & 0 & \cdots & 0 \\
  \mu &  0  & 0 & \cdots & 0 \\
  0   &  0  & 0 & \cdots & 0 \\
  \vdots & \vdots & \vdots & \ddots & \vdots \\
  0 & 0 & 0 & \cdots & 0\\
  \end{pmatrix}  , \   
 A_2  =  \begin{pmatrix}
    \mu & 0 & 0 & \cdots & 0 \\
  0 &   -\mu & 0 & \cdots & 0 \\
  0 & 0 & 0 & \cdots & 0 \\
  \vdots & \vdots & \vdots & \ddots & \vdots \\
  0 & 0 & 0 & \cdots & 0 \\
  \end{pmatrix}
, \  
 A_3  =  \begin{pmatrix}
  c &     0  & 0 & \cdots & 0 \\
  0      & c & 0 & \cdots & 0 \\
  0 & 0 & c & \cdots & 0 \\
  \vdots & \vdots & \vdots & \ddots & \vdots \\
  0 & 0 & 0 & \cdots & c \\
  \end{pmatrix}
,     \\
A_4  &= \ldots = A_m = 0 ,
\end{aligned}
 \end{equation}
 where
$ c $ and
 $  \mu $ are real functions on $ M $ such that $ \mu \ne 0 $
and $ c^2 = - {\widetilde k} $.
Moreover, $M$ is Weyl-semi-symmetric
in $ \R^{n+m}({\widetilde k})$~; and,
in   the second   case,
$M$ is a minimal or
 pseudo-umbilical submanifold
  in $ \R^{n+m}({\widetilde k})$. \   $\square$
       \end{theoreme}

   \begin{corollaire}   \label{cr.LQgCorigin618841}   
  Let $ M $  be  a Wintgen ideal
	submanifold of codimension
  $ m $ in a real space form    $   \R^{n+m}({\widetilde k})  $,
  $ n \geq 4 $ and  $ m \geq 1 $.
	  If $ {\widetilde k} = 0$,
	then  the tensors
	$ R \cdot C $ and Q(\mbox{\rm Ricc},C)
	 are linearly dependent
		if and only if (i) either $ M$  is totally umbilical
in  $   \R^{n+m}  $
(and hence  $M$ is conformally flat),
(ii) or     with respect to some suitable adapted
Choi-Lu frames
  $  \left\{E_i , \xi_\alpha\right\}    $
  on $ M $ in $ \R^{n+m} $, the
  shape operators are given by
\begin{equation}\label{ChoiLu6194112}
\begin{aligned}
  A_1 &=  \begin{pmatrix}
  0   & \mu & 0 & \cdots & 0 \\
  \mu &  0  & 0 & \cdots & 0 \\
  0   &  0  & 0 & \cdots & 0 \\
  \vdots & \vdots & \vdots & \ddots & \vdots \\
  0 & 0 & 0 & \cdots & 0\\
  \end{pmatrix}  , \   
 A_2  =  \begin{pmatrix}
    \mu & 0 & 0 & \cdots & 0 \\
  0 &   -\mu & 0 & \cdots & 0 \\
  0 & 0 & 0 & \cdots & 0 \\
  \vdots & \vdots & \vdots & \ddots & \vdots \\
  0 & 0 & 0 & \cdots & 0 \\
  \end{pmatrix}
, \\
 A_3 &=  A_4  = \ldots = A_m = 0 ,
\end{aligned}
 \end{equation}
 where
 $\mu $ is a real function
 on $ M $ such that $ \mu \ne 0 $.
  In this second case,  $M$ is minimal
  in $ \R^{n+m} $. \ $\square$
  \  $\square$

      \end{corollaire}

  \begin{theoreme}   \label{cr.LQgCorigin61882}   
  Let $ M = M^n$  be  a Wintgen ideal
	submanifold of codimension
  $ m $ in a real space form
	$ \R^{n+m}({\widetilde k})  $,
  $ n \geq 4 $ and  $ m \geq 1 $.
	Then   the tensors
	$ C \cdot R $ and $Q(\mbox{\rm Ricc}, C)$
 are linearly dependent	if and only if
	$ M$  is totally umbilical
in  $   \R^{n+m}({\widetilde k})   $
(and hence  $M$ is conformally flat).
 \    $  \square $
       \end{theoreme}

  \begin{theoreme}   \label{cr.LQgCorigin619819}   
  Let $ M = M^n $  be  a Wintgen ideal
	submanifold of codimension
  $ m $ in a real space form
	$ \R^{n+m}({\widetilde k}) $,
  $ n \geq 4 $ and  $ m \geq 1 $.
	Then   the tensors
	$ R \cdot C  - C \cdot R$ and $Q(\mbox{\rm Ricc} , C)$
	 are linearly dependent
if and only if
  $ M$  is totally umbilical
in  $   \R^{n+m}({\widetilde k})  $
(and hence  $M$ is conformally flat). \      $\square$
        \end{theoreme}




   \begin{theoreme}   \label{cr.LQgCorigin61844110}   
  Let $ M = M^n $  be  a Wintgen ideal
	submanifold of codimension
  $ m $ in a real space form
	$  \R^{n+m}({\widetilde k})  $,
  $ n \geq 4 $ and  $ m \geq 1 $.
Then   the tensors
	$ R \cdot C $ and $Q(\mbox{\rm Ricc} , g\wedge\mbox{\rm Ricc})$
		 are linearly dependent
			if and only if (i) either $ M$  is totally umbilical
in  $   \R^{n+m}({\widetilde k})   $
(and hence  $M$ is conformally flat),
(ii) or   with respect to some suitable adapted
Choi-Lu frames
  $  \left\{E_i , \xi_\alpha\right\}    $
  on $ M $ in $ \R^{n+m}({\widetilde k})  $, the
  shape operators are given by
\begin{equation}\label{ChoiLu618413}
\begin{aligned}
  A_1 &=  \begin{pmatrix}
  0   & \mu & 0 & \cdots & 0 \\
  \mu &  0  & 0 & \cdots & 0 \\
  0   &  0  & 0 & \cdots & 0 \\
  \vdots & \vdots & \vdots & \ddots & \vdots \\
  0 & 0 & 0 & \cdots & 0\\
  \end{pmatrix}  , \    
 A_2  =  \begin{pmatrix}
    \mu & 0 & 0 & \cdots & 0 \\
  0 &   -\mu & 0 & \cdots & 0 \\
  0 & 0 & 0 & \cdots & 0 \\
  \vdots & \vdots & \vdots & \ddots & \vdots \\
  0 & 0 & 0 & \cdots & 0 \\
  \end{pmatrix}
, \   
 A_3  =  \begin{pmatrix}
  c &     0  & 0 & \cdots & 0 \\
  0      & c & 0 & \cdots & 0 \\
  0 & 0 & c & \cdots & 0 \\
  \vdots & \vdots & \vdots & \ddots & \vdots \\
  0 & 0 & 0 & \cdots & c \\
  \end{pmatrix}
,    \\
A_4  &= \ldots = A_m = 0 ,
\end{aligned}
 \end{equation}
 where
$ c $ and
 $  \mu $ are real functions on $ M $ such that $ \mu \ne 0 $
 and $   c^2 + {\widetilde k} \ne 0 $.
Moreover,
 $$
   R \cdot C
   =  \frac{{\widetilde k} + c^2}{2\mu^2}  Q(\mbox{\rm Ricc} ,
   g\wedge\mbox{\rm Ricc}) ~;
 $$
and,
in the second  case,
 $M$ is a minimal
 or    pseudo-umbilical submanifold
in $ \R^{n+m}({\widetilde k})$. \  $\square$
       \end{theoreme}

  \begin{theoreme}   \label{cr.LQgCorigin618441}   
  Let $ M = M^n$  be  a Wintgen ideal
	submanifold of codimension
  $ m $ in a real space form
	$ \R^{n+m}({\widetilde k})  $,
  $ n \geq 4 $ and  $ m \geq 1 $.
	Then   the tensors
	$ C \cdot R $ and $Q(\mbox{\rm Ricc}, g\wedge\mbox{\rm Ricc})$
 are linearly dependent	if and only if
	$ M$  is totally umbilical
in  $   \R^{n+m}({\widetilde k})   $
(and hence  $M$ is conformally flat).
 \    $  \square $
       \end{theoreme}

  \begin{theoreme}   \label{cr.LQgCorigin614411}   
  Let $ M = M^n $  be  a Wintgen ideal
	submanifold of codimension
  $ m $ in a real space form
	$ \R^{n+m}({\widetilde k})  $,
  $ n \geq 4 $ and  $ m \geq 1 $.
	Then   the tensors
	$ R \cdot C  - C \cdot R$ and $Q(\mbox{\rm Ricc} , g\wedge\mbox{\rm Ricc})$
	 are linearly dependent
if and only if
(i) either $ M$  is totally umbilical
in  $   \R^{n+m}({\widetilde k})  $
(and hence  $M$ is conformally flat),
(ii) or   with respect to some suitable adapted
Choi-Lu frames
  $  \left\{E_i , \xi_\alpha\right\}    $
  on $ M $ in $ \R^{n+m}({\widetilde k}) $, the
  shape operators are given by
\begin{equation}\label{ChoiLu61841212}
\begin{aligned}
  A_1 &=  \begin{pmatrix}
  0   & \mu & 0 & \cdots & 0 \\
  \mu &  0  & 0 & \cdots & 0 \\
  0   &  0  & 0 & \cdots & 0 \\
  \vdots & \vdots & \vdots & \ddots & \vdots \\
  0 & 0 & 0 & \cdots & 0\\
  \end{pmatrix}  , \    
 A_2  =  \begin{pmatrix}
     2\mu          &    0         &      0     & \cdots & 0 \\
     0             &    0         &      0     & \cdots & 0 \\
     0             &    0         &    \mu     &    0   &  0 \\
  \vdots           &   \vdots     &    \vdots  & \ddots & \vdots \\
     0             &    0         &      0     & \cdots & \mu \\
  \end{pmatrix}
, \    
  A_3  =  \begin{pmatrix}
  c &     0  & 0 & \cdots & 0 \\
  0      & c & 0 & \cdots & 0 \\
  0 & 0 & c & \cdots & 0 \\
  \vdots & \vdots & \vdots & \ddots & \vdots \\
  0 & 0 & 0 & \cdots & c \\
  \end{pmatrix}
,     \\
 A_4   &= \cdots = A_m = 0 ,
\end{aligned}
 \end{equation}
 where $c$ and
 $  \mu $ are real functions  on $ M $ such that $ \mu \ne 0 $~;
in this  second case,
$$
H^2 + {\widetilde k} = \frac{2}{n-1}\mu^2
    \  \ \text{and} \ \
R \cdot C  - C \cdot R =
 -\frac{n-3}{(n-1)(n-2)} Q(\mbox{\rm Ricc} , g\wedge\mbox{\rm Ricc})
 . \      \square
 $$
       \end{theoreme}

  \begin{corollaire}   
  Let $ M = M^n$  be
 a  Wintgen ideal   submanifold
  of codimension $ m $
   in a real space form
	$  \R^{n+m}({\widetilde k})  $,
	$ n \geq 4 $ and  $ m \geq 1 $. Suppose $ {\widetilde k} \leq 0$.
	The following  
	 conditions  are both equivalent.
\begin{enumerate}
 \item[(I)]
   With respect to some suitable adapted
Choi-Lu frames
  $  \left\{E_i , \xi_\alpha\right\}    $
  on $ M $ in $ \R^{n+m}({\widetilde k})  $, the
  shape operators are given by
\begin{equation}\label{ChoiLu231}
\begin{aligned}
  A_1 &=  \begin{pmatrix}
  0   & \mu & 0 & \cdots & 0 \\
  \mu &  0  & 0 & \cdots & 0 \\
  0   &  0  & 0 & \cdots & 0 \\
  \vdots & \vdots & \vdots & \ddots & \vdots \\
  0 & 0 & 0 & \cdots & 0\\
  \end{pmatrix}  , \
 A_2  =  \begin{pmatrix}
    \mu & 0 & 0 & \cdots & 0 \\
  0 &   -\mu & 0 & \cdots & 0 \\
  0 & 0 & 0 & \cdots & 0 \\
  \vdots & \vdots & \vdots & \ddots & \vdots \\
  0 & 0 & 0 & \cdots & 0 \\
  \end{pmatrix}
, \   
 A_3   =  \begin{pmatrix}
  c &     0  & 0 & \cdots & 0 \\
  0      & c & 0 & \cdots & 0 \\
  0 & 0 & c & \cdots & 0 \\
  \vdots & \vdots & \vdots & \ddots & \vdots \\
  0 & 0 & 0 & \cdots & c \\
  \end{pmatrix}
,   \\
A_4   &= \ldots = A_m = 0 ,
\end{aligned}
 \end{equation}
 where
$ c $ and
 $  \mu $ are real functions on $ M $ such that $ \mu \ne 0 $
 and $ c^2 = - {\widetilde k} $.

\item[(II)]     $   M   $ is Weyl-pseudo-symmetric, i.e.,
$  R \cdot C   $   and $ Q(g, R)  $ are linearly  dependent.

\item[(III)]   $  R \cdot C  $
and $ Q(g, g\wedge \mbox{\rm Ricc})  $ are linearly  dependent.

	\item[(IV)]    $  R \cdot C  $
and $ Q(\mbox{\rm Ricc} , R)  $ are linearly  dependent.

\item[(V)]       $  R \cdot C  $
and $ Q(\mbox{\rm Ricc} , C)  $ are linearly  dependent.

\item[(VI)]  $  R \cdot C  - C \cdot R $
and $ Q(g, R)  $ are linearly  dependent. 

\end{enumerate}

\noindent In these cases (I) to  (V),
$M$ is Weyl-semi-symmetric
in $ \R^{n+m}({\widetilde k})  $.
In the case (VI),
$$
   R \cdot C - C \cdot R
   =  -\dfrac{2(n-3)\mu^2}{(n-1)(n-2)} Q(g, R)  . \    \square
   $$

   \end{corollaire}


  \begin{corollaire}   
  Let $ M = M^n$  be
 a  Wintgen ideal   submanifold
  of codimension $ m $
   in a real space form
	$  \R^{n+m}({\widetilde k}) $,
	$ n \geq 4 $ and  $ m \geq 1 $.
	Suppose $ {\widetilde k} = 0$.
	The following  
	 conditions  are both equivalent.
\begin{enumerate}
 \item[(I)]
  $M$ is minimal and  with respect to some suitable adapted
Choi-Lu frames
  $  \left\{E_i , \xi_\alpha\right\}    $
  on $ M $ in $ \R^{n+m}({\widetilde k})  $, the
  shape operators are given by
\begin{equation}\label{ChoiLu23190}
\begin{aligned}
  A_1 &=  \begin{pmatrix}
  0   & \mu & 0 & \cdots & 0 \\
  \mu &  0  & 0 & \cdots & 0 \\
  0   &  0  & 0 & \cdots & 0 \\
  \vdots & \vdots & \vdots & \ddots & \vdots \\
  0 & 0 & 0 & \cdots & 0\\
  \end{pmatrix}  , \    
 A_2  =  \begin{pmatrix}
    \mu & 0 & 0 & \cdots & 0 \\
  0 &   -\mu & 0 & \cdots & 0 \\
  0 & 0 & 0 & \cdots & 0 \\
  \vdots & \vdots & \vdots & \ddots & \vdots \\
  0 & 0 & 0 & \cdots & 0 \\
  \end{pmatrix}
, \\
 A_3 &= A_4  = \ldots = A_m = 0 ,
\end{aligned}
 \end{equation}
 where
  $  \mu $ is a real function  on $ M $ such that $ \mu \ne 0 $.

\item[(II)]     $   M   $ is Weyl-pseudo-symmetric, i.e.
$  R \cdot C   $   and $ Q(g, R)  $ are linearly  dependent.

\item[(III)]   $  R \cdot C  $
and $ Q(g, g\wedge \mbox{\rm Ricc})  $ are linearly  dependent.

	\item[(IV)]    $  R \cdot C  $
and $ Q(\mbox{\rm Ricc} , R)  $ are linearly  dependent.

\item[(V)]       $  R \cdot C  $
and $ Q(\mbox{\rm Ricc} , C)  $ are linearly  dependent.

\item[(VI)]  $  R \cdot C  - C \cdot R $
and $ Q(g, R)  $ are linearly  dependent.  \      $\square$

\end{enumerate}

   \end{corollaire}



  \begin{corollaire}   
  Let $ M = M^n$  be
 a  Wintgen ideal   submanifold
  of codimension $ m $
   in a real space form
	$  \R^{n+m}({\widetilde k})  $,
	$ n \geq 4 $ and  $ m \geq 1 $.
	Suppose $ {\widetilde k} > 0$.
	Then the following conditions are both equivalent.
\begin{enumerate}
 \item[(I)] $  R \cdot C   $   and $ Q(g, R)  $ are   linearly  dependent.

\item[(II)]   $  R \cdot C   $
and $ Q(g, g\wedge \mbox{\rm Ricc})  $ are   linearly  dependent.

\item[(III)]  $  R \cdot C  $  and $ Q(\mbox{\rm Ricc}, R)  $
are   linearly  dependent.

\item[(IV)]    $  R \cdot C  $  and $ Q(\mbox{\rm Ricc}, C)  $
are   linearly  dependent.

\item[(V)]    $  M  $  is totally umbilical in $  \R^{n+m}({\widetilde k})  $
(and hence $M$ is conformally flat).   \      $\square$

\end{enumerate}

   \end{corollaire}

\subsection{Proofs of main results}

\subsubsection{Introduction}

We consider  mainly the Wintgen ideal 
submanifold
$ M = M^n $ of dimension $ n > 3 $
and of codimension $ m \geq 2 $ in a real space form
$  \R^{n+m}({\widetilde k})  $
 and the  Theorem $A$ in the last section  \ref{ciso2}.
	We note that
	the second fundamental tensor $ H  $
   satisfies:
	 $$
	H^2 =	a^2 + b^2 + c^2  .
	  $$
		As in {\bf Theorem $A$},
 the  shape operator $ A $ satifies the following relations,
  for all $ \beta \geq 4 $
 and $ 3 \leq i < j \leq n $~:
\begin{eqnarray}
    &  &
    \left\{
  \begin{aligned}
   A_1(E_1) \wedge_g A_1(E_2)  &= \left(a^2 - \mu^2\right)
    E_1 \wedge_g  E_2  , \\
   A_2(E_1) \wedge_g A_2(E_2)  &= (b^2 - \mu^2)E_1 \wedge_g  E_2  ,  \\
    A_3(E_1) \wedge_g A_3(E_2)  &= c^2E_1 \wedge_g  E_2   , \\
         A_\beta(E_1) \wedge_g A_\beta(E_2) &= 0   , \\
  \end{aligned}
     \right.
 \label{eqbli3}  
 \end{eqnarray}
 \begin{eqnarray}
   &  &
   \left\{
  \begin{aligned}
   A_1(E_1) \wedge_g A_1(E_i)  &= a^2 E_1 \wedge_g  E_i +  a\mu E_2 \wedge_g  E_i  , \\
      A_2(E_1) \wedge_g A_2(E_i)  &= (b^2 + b\mu)E_1 \wedge_g  E_i   , \\
   A_3(E_1) \wedge_g A_3(E_i) &= c^2E_1 \wedge_g  E_i   , \\
         A_\beta(E_1) \wedge_g A_\beta(E_i) &= 0   , \\
  \end{aligned}
     \right.
 \label{eqbli4} 
 \end{eqnarray}
 \begin{eqnarray}
   &  &
 \left\{
  \begin{aligned}
    A_1(E_2) \wedge_g A_1(E_i)  &=   a\mu E_1 \wedge_g  E_i +  a^2 E_2 \wedge_g  E_i  , \\
      A_2(E_2) \wedge_g A_2(E_i)  &= (b^2 - b\mu)E_2 \wedge_g  E_i   , \\
   A_3(E_2) \wedge_g A_3(E_i) &= c^2E_2 \wedge_g  E_i  , \\
        A_\beta(E_2) \wedge_g A_\beta(E_i) &= 0   , \\
   \end{aligned}
     \right.
\label{eqbli5} 
\end{eqnarray}
 \begin{eqnarray}
   &  &
 \left\{
  \begin{aligned}
    A_1(E_i) \wedge_g A_1(E_j)  &=    a^2 E_i \wedge_g  E_j  , \\
      A_2(E_i) \wedge_g A_2(E_j)  &= b^2 E_i \wedge_g  E_j   , \\
   A_3(E_i) \wedge_g A_3(E_j) &= c^2 E_i \wedge_g  E_j   , \\
      A_\beta(E_i) \wedge_g A_\beta(E_j) &= 0   . \\
  \end{aligned}
     \right.
     \label{eqbli6}
  \end{eqnarray}

Due to the relations (\ref{eqbli3}) to  (\ref{eqbli6}),
  the Riemannian-Christoffel curvature operator $ \mathcal{R} $
  satisfies the following equalities~:
 \begin{equation}\label{eq010}
 \left\{
  \begin{aligned}
 \mathcal{R}(E_1 , E_2) &= \left(H^2 + {\widetilde k} - 2\mu^2\right) E_1\wedge_g E_2  , \\
\mathcal{R}(E_1 , E_i) &= \left(H^2
              + {\widetilde k} + b\mu\right)E_1\wedge_g E_i + a\mu E_2\wedge_g E_i  , \\
 \mathcal{R}(E_2 , E_i) &= a\mu E_1\wedge_g E_i
              + \left(H^2  + {\widetilde k} - b\mu\right) E_2\wedge_g E_i  , \\
  \mathcal{R}(E_i , E_j) &= (H^2 + {\widetilde k})E_i\wedge_g E_j   . \\
     \end{aligned}
     \right.
 \end{equation}

\noindent
 The local components
$ R_{u v  w  t} =  R\left(E_u , E_v , E_w , E_t\right)  $
 of the Riemannian-Christoffel curvature $ (0,4) $-tensor $ R $
     are given by~:
 \begin{equation}\label{eq001*}
 \left\{
  \begin{aligned}
  R_{1221}     &=   H^2   + {\widetilde k}  - 2\mu^2   , \\
  R_{1ii1}     &=   H^2  + {\widetilde k} + b\mu \quad \text{\rm for} \quad i \geq 3 ,  \\
  R_{2ii2}     &=  H^2  + {\widetilde k}  - b\mu  \quad \text{\rm for} \quad i \geq 3 , \\
  R_{ijji}     &=    H^2   +  {\widetilde k} \quad \text{\rm for} \quad 3 \leq i < j \leq n  , \\
   R_{1ii2}     &=      a\mu    \quad \text{\rm for} \quad 3 \leq i   \leq n , \\
    \end{aligned}
     \right.
 \end{equation}
the other values of $  R_{u v  w  t}  $ being null.

    We denote $ \mbox{\rm Ricc} $   the Ricci tensor.
     We set $S_{u v} = \mbox{\rm Ricc}\left(E_u , E_v\right)$. 
Then~:
   \begin{equation}\label{eq001**}
 \left\{
  \begin{aligned}
  S_{11}  &=    (n-1) \left(H^2   + {\widetilde k}\right)  - 2\mu^2 + (n-2)b\mu   , \\
  S_{12}  &=    (n-2)a\mu     , \\
  S_{22}  &= (n-1)\left(H^2   + {\widetilde k}\right)    - 2\mu^2 -(n-2)b\mu
	                                      \quad \text{\rm for} \quad i \geq 3 ,  \\
  S_{ii}  &=   (n-1)\left(H^2   + {\widetilde k}\right)
	                                  \quad \text{\rm for} \quad 3 \leq i \leq n    ,  \\
     \end{aligned}
     \right.
 \end{equation}
the other values of $S_{u v}$ being null. Next, computing the Ricci operator 
$ \mathcal{S} $  (associated  to the Ricci tensor $ \mbox{\rm Ricc} $) 
by using the above equalities, we obtain~:
 \begin{equation}\label{eq3}
 \left\{
  \begin{aligned}
\mathcal{S}(E_1) &= \mbox{\rm Ricc}\left(E_1,E_1\right)E_1 + \mbox{\rm Ricc}\left(E_1,E_2\right)E_2 \\
&= \left[(n-1)\left(H^2 + {\widetilde k}\right) - 2\mu^2 + (n-2)b\mu\right]E_1
               +  \left[(n-2)a\mu\right]E_2 , \\
\mathcal{S}(E_2) &= \mbox{\rm Ricc}\left(E_2,E_1\right)E_1 + \mbox{\rm Ricc}\left(E_2,E_2\right)E_2 \\
&=    \left[(n-2)a\mu\right]E_1
  +  \left[ (n-1)\left(H^2   + {\widetilde k}\right)   - 2\mu^2 -  (n-2)b\mu\right]E_2  , \\
\mathcal{S}(E_i) &= \mbox{\rm Ricc}\left(E_i,E_1\right)E_1 + \mbox{\rm Ricc}\left(E_i,E_2\right)E_2 \\
&=   \left[(n-1)\left(H^2   + {\widetilde k}\right) \right]E_i
			 \quad \text{\rm for} \quad i \geq 3 .  \\
      \end{aligned}
      \right.
 \end{equation}

The scalar curvature $ \tau $ is given by~:
 \begin{equation}\label{eq5}
\tau = \sum_{i=1}^n \mbox{\rm Ricc}(E_i , E_i)
 = n(n-1) \left(H^2   + {\widetilde k}\right)  - 4\mu^2 ,
\end{equation}
 so that the normalized scalar of $M$ is
  \begin{equation}\label{eq6}
 \rho = \dfrac{\tau}{n(n-1)} = \left(H^2   + {\widetilde k}\right)  -  \dfrac{4\mu^2}{n(n-1)}  .
\end{equation}

We recall that for any tangent vector fields $X$, $Y$, $Z$, $W$ we have
$$
\begin{aligned}
	(g\wedge\mbox{\rm Ricc})(X ,  Y ;  Z ,  W)
&   = g(X , W)\mbox{\rm Ricc}(Y,Z) + \mbox{\rm Ricc}(X,W)g(Y,Z) \\
&  -  g(X , Z)\mbox{\rm Ricc}(Y,W) - \mbox{\rm Ricc}(X,Z)g(Y,W) \\
&= g\left(\left(X\wedge_g  \mathcal{S}(Y) +  \mathcal{S}(X)\wedge_g Y\right)(Z) ,   W\right) .\\
 \end{aligned}
$$

\noindent If  we set~:
$  \left(g\wedge\mbox{\rm Ricc}\right)_{u  v  s  t} =
  \left(g\wedge\mbox{\rm Ricc}\right)\left(E_u , E_v , E_s , E_t\right)$,
for any $u \geq 1$, $v \geq 1$, $s \geq 1$ and $t \geq 1$,
then~:
   \begin{equation}\label{eqgWRICCI001**}
 \left\{
  \begin{aligned}
\left(g\wedge\mbox{\rm Ricc}\right)_{1221} &= 2(n-1)\left(H^2 + {\widetilde k}\right) - 4\mu^2 , \\
\left(g\wedge\mbox{\rm Ricc}\right)_{1ii1} &= 2(n-1) \left(H^2 + {\widetilde k}\right)
                                                                       - 2\mu^2 + (n-2)b\mu , \\
  \left(g\wedge\mbox{\rm Ricc}\right)_{2ii2}    &= 2(n-1)({\widetilde k}  + H^2) - 2\mu^2 - (n-2)b\mu
		                                                       \quad \text{\rm for} \quad i \geq 3 ,  \\
  \left(g\wedge\mbox{\rm Ricc}\right)_{ijji}   &=  2(n-1)\left(H^2   + {\widetilde k}\right)
		                                              \quad \text{\rm for} \quad 3 \leq i \leq n    ,  \\
     \end{aligned}
     \right.
 \end{equation}
the other values of $  \left(g\wedge\mbox{\rm Ricc}\right)_{u v st} $ being null.

Let $ C $ be the Weyl conformal curvature tensor.
Putting $C_{u v  w  t} = C\left(E_u , E_v , E_w , E_t\right)$, we get~:
 \begin{equation}\label{eq701*}
 \left\{
  \begin{aligned}
  C_{1221}     &= - \dfrac{2(n - 3)\mu^2}{n - 1}    , \\
  C_{1ii1}     &=   \dfrac{2(n - 3)\mu^2}{(n - 1)(n-2)}    \quad \text{\rm for} \quad i \geq 3 ,  \\
  C_{2ii2}     &=    \dfrac{2(n - 3)\mu^2}{(n - 1)(n-2)}   \quad \text{\rm for} \quad i \geq 3 , \\
  C_{ijji}     &=   - \dfrac{4\mu^2}{(n - 1)(n-2)}   \quad \text{\rm for} \quad 3 \leq i < j \leq n ,  \\
     \end{aligned}
     \right.
 \end{equation}
the other values of $  C_{u v  w  t}  $ being null.

\subsubsection{Proofs of
theorems  {\sl \bf 1, 2, 3, 4}
and  corollary    {\sl \bf 1}}

Now we compute the local components of the tensors 
   $ R \cdot C  $,
	 $   C \cdot R $ and
	$ R \cdot C - C \cdot R $ of a   Wintgen  ideal
  submanifold.
	  Let  
		$Z$, $W$ be  tangent  vector fields.

Firstly, we compute the local components of the tensor  $ R \cdot C$.
For   any index $ i \in \left\{3, \ldots , n\right\}  $,
  \begin{equation}  \label{r.c_01}
   \left\{
\begin{aligned}
&
\left(R \cdot C\right)\left(E_1 , E_2 , Z , W; E_1 , E_i\right)
=   \dfrac{2(n-3)\mu^3a}{n-2}  
 \left\langle \left(E_1 \wedge_g E_i\right)(Z) , W\right\rangle  \\
&
 - \dfrac{2(n-3)\mu^2}{n-2}\left(H^2  + {\widetilde k} + \mu b\right)
					\left\langle \left(E_2 \wedge_g E_i\right)(Z) , W\right\rangle , \\
&
\left(R \cdot C\right)\left(E_1 , E_2 , Z , W; E_2 , E_i\right)
=   \dfrac{2(n-3)\mu^2}{n-2}\left(H^2 + {\widetilde k} - \mu b\right)
               \left\langle \left(E_1 \wedge_g E_i\right)(Z) , W\right\rangle \\
 &
  - \dfrac{2(n-3)\mu^3a}{n-2}\left\langle \left(E_2 \wedge_g E_i\right)(Z) , 
  W\right\rangle  . \\
\end{aligned}
 	\right.
\end{equation}
\noindent For   any   indexes
   $ i  $, $j   \in \left\{3 , \ldots ,  n\right\}$
	such that $ i \not= j$,  	
  \begin{equation} \label{r.c_02}
 \left\{
  \begin{aligned}
 \left(R \cdot C\right)\left(E_1 , E_i , Z , W; E_1 , E_2\right)
 &=   0  , \\
 \left(R \cdot C\right)\left(E_1 , E_i , Z , W; E_1 , E_j\right)
 &=    \dfrac{2\mu^2}{n-2}\left(H^2 + {\widetilde k}
                    + \mu b\right)
                    \left\langle \left(E_i \wedge_g E_j\right)(Z) , W\right\rangle , \\
 \left(R \cdot C\right)\left(E_1 , E_i , Z , W; E_2 , E_i\right)
 &=   \dfrac{2(n-3)\mu^2}{n-2}
 \left(H^2 + {\widetilde k} - \mu b\right)
 \left\langle \left(E_1 \wedge_g E_2\right)(Z) , W\right\rangle , \\
 \left(R \cdot C\right)\left(E_1 , E_i , Z , W; E_i , E_j\right)
 &=    0 . \\
  \end{aligned}
	\right.
 \end{equation}

    \noindent   For   any    indexes
   $ i  $, $j  \in \left\{3 , \ldots ,  n\right\}$
	such that $ i \not= j$,
  \begin{equation} \label{r.c_03}
 \left\{
  \begin{aligned}
  \left(R \cdot C\right)\left(E_2 , E_i , Z , W; E_2 , E_1\right)
 &=   0  , \\
 \left(R \cdot C\right)\left(E_2 , E_i , Z , W; E_2 , E_j\right)
 &=    \dfrac{2\mu^2}{n-2}\left(H^2 + {\widetilde k} - \mu b\right)
\left\langle \left(E_i \wedge_g E_j\right)(Z) , W\right\rangle , \\
 \left(R \cdot C\right)\left(E_2 , E_i , Z , W; E_i , E_1\right)
 &=    \dfrac{2(n-3)\mu^2}{n-2}\left[H^2 + {\widetilde k}   + \mu b\right]
\left\langle \left(E_1 \wedge_g E_2\right)(Z) , W\right\rangle  , \\
\left(R \cdot C\right)\left(E_2 , E_i , Z , W; E_i , E_j\right)
 &=   0 .  \\
  \end{aligned}
	\right.
 \end{equation}

    \noindent   For   any     indexes
   $ i  $, $j$,  $ k    \in \left\{3 , \ldots ,  n\right\}$
	such that $ i \not= j$, $ i \not= k$ and $ j \not= k$,
  \begin{equation}  \label{r.c_04}
 \left\{
  \begin{aligned}
 \left(R \cdot C\right)\left(E_i , E_j , Z , W; E_i , E_1\right)
 &=    \dfrac{2\mu^2}{n-2}\left(H^2 + {\widetilde k} + \mu b\right)\left\langle
 \left(E_1 \wedge_g E_j\right)(Z) , W\right\rangle \\
& \hspace{0,2cm}
+
\dfrac{2\mu^2a}{n-2}\left\langle\left(E_2 \wedge_g E_j\right)(Z) , W\right\rangle , \\
\left(R \cdot C\right)\left(E_i , E_j , Z , W; E_i , E_2\right)
 &= \dfrac{2\mu^3a}{n-2}\left\langle\left(E_1 \wedge_g E_j\right)(Z) , W\right\rangle \\
& \hspace{0,2cm}
+
\dfrac{2\mu^2}{n-2}\left(H^2 + {\widetilde k} - \mu b\right)
\left\langle\left(E_2 \wedge_g E_j\right)(Z) , W\right\rangle , \\
\left(R \cdot C\right)\left(E_i , E_j , Z , W; E_i , E_k\right)
 &=    0 , \\
\left(R \cdot C\right)\left(E_i , E_j , Z , W; E_j , E_k\right)
 &=    0  . \\
  \end{aligned}
	\right.
 \end{equation}

Secondly, we compute the local components of the tensor $ C \cdot R$.
For any index
     $ i   \in \left\{3 , \ldots ,  n\right\}$,
  \begin{equation}  \label{c.r_01}
 \left\{
  \begin{aligned}
\left(C \cdot R\right)\left(E_1 , E_2 , Z , W; E_1 , E_i\right)
 &=    \dfrac{2(n-3)\mu^3a}{(n-1)(n-2)}\left\langle \left(E_1 \wedge_g E_i\right)(Z), W\right\rangle \\
& \hspace{0,2cm} +
 \dfrac{2(n-3)\mu^2}{(n-1)(n-2)}\left(2\mu^2 - \mu b\right)
                    \left\langle \left(E_2 \wedge_g E_i\right)(Z) , W\right\rangle , \\
\left(C \cdot R\right)\left(E_1 , E_2 , Z , W; E_2 , E_i\right)
 &=   \dfrac{2(n-3)\mu^2}{(n-1)(n-2)}\left(2\mu^2 + \mu b\right)
                     \left\langle \left(E_1 \wedge_g E_i\right)(Z) , W\right\rangle \\
& \hspace{0,2cm}  +
 \dfrac{2(n-3)\mu^3a}{(n-1)(n-2)}
                     \left\langle \left(E_2 \wedge_g E_i\right)(Z) , W\right\rangle . \\
   \end{aligned}
	\right.
    \end{equation}

\noindent   For   any    indexes
   $ i  $, $j    \in \left\{3 , \ldots ,  n\right\}$
	such that $ i \not= j$,
  \begin{equation}  \label{c.r_02}
 \left\{
  \begin{aligned}
\left(C \cdot R\right)\left(E_1 , E_i , Z , W; E_1 , E_2\right)
 &=   -\dfrac{4(n-3)\mu^3a}{n-1}  \left\langle \left(E_1 \wedge_g E_i\right)(Z) , W\right\rangle \\
 & \hspace{0,2cm}
  +
 \dfrac{4(n-3)\mu^3b}{n-1}\left\langle \left(E_2 \wedge_g E_i\right)(Z) , W\right\rangle , \\
\left(C \cdot R\right)\left(E_1 , E_i , Z , W; E_1 , E_j\right)
 &=    \dfrac{2(n-3)\mu^3b}{(n-1)(n-2)}\left\langle \left(E_i \wedge_g E_j\right)(Z) , W\right\rangle , \\
 \left(C \cdot R\right)\left(E_1 , E_i , Z , W; E_2 , E_i\right)
 &=     \dfrac{2(n-3)\mu^2}{(n-1)(n-2)}\left(2\mu^2 + \mu b\right)\left\langle \left(E_1 \wedge_g E_2\right)(Z) , W\right\rangle , \\
 \left(C \cdot R\right)\left(E_1 , E_i , Z , W; E_i , E_j\right)
 &=    0  . \\
   \end{aligned}
	\right.
    \end{equation}

\noindent   For   any indexes
   $ i  $, $j   \in \left\{3 , \ldots ,  n\right\}$
	such that $ i \not= j$,
 \begin{equation}  \label{c.r_03}
 \left\{
  \begin{aligned}
 \left(C \cdot R\right)\left(E_2 , E_i , Z , W; E_2 , E_1\right)
 &=    -\dfrac{4(n-3)\mu^3b}{n-1}\left\langle \left(E_1 \wedge_g E_i\right)(Z) , W\right\rangle \\
& \hspace{0,2cm} 
 -
 \dfrac{4(n-3)\mu^3a}{n-1}\left\langle \left(E_2 \wedge_g E_i\right)(Z) , W\right\rangle , \\
 \left(C \cdot R\right)\left(E_2 , E_i , Z , W; E_2 , E_j\right)
 &=	-\dfrac{2(n-3)\mu^3b}{(n-1)(n-2)}\left\langle \left(E_i \wedge_g E_j\right)(Z) , W\right\rangle , \\
 \left(C \cdot R\right)\left(E_2 , E_i , Z , W; E_i , E_1\right)
 &=	\dfrac{2(n-3)\mu^2}{(n-1)(n-2)}\left(2\mu^2 - b\mu\right)
\left\langle \left(E_1 \wedge_g E_2\right)(Z) , W\right\rangle , \\
\left(C \cdot R\right)\left(E_2 , E_i , Z , W; E_i , E_j\right)
 &=	0 . \\
   \end{aligned}
	\right.
    \end{equation}

\noindent   For   any indexes
   $ i  $, $j$,  $ k    \in \left\{3 , \ldots ,  n\right\}$
	such that $ i \not= j$, $ i \not= k$ and $ j \not= k$,
\begin{equation}  \label{c.r_04}
 \left\{
  \begin{aligned}
 \left(C \cdot R\right)\left(E_i , E_j , Z , W; E_i , E_1\right)
 &=    \dfrac{2(n-3)\mu^3b}{(n-1)(n-2)}\left\langle \left(E_1 \wedge_g E_j\right)(Z) , W\right\rangle \\
& \hspace{0,2cm}  +
	\dfrac{2(n-3)\mu^3a}{(n-1)(n-2)}\left\langle \left(E_2 \wedge_g E_j\right)(Z) , W\right\rangle , \\
 \left(C \cdot R\right)\left(E_i , E_j , Z , W; E_i , E_2\right)
 &=     \dfrac{2(n-3)\mu^3a}{(n-1)(n-2)}\left\langle \left(E_1 \wedge_g E_i\right)(Z) , W\right\rangle \\
& \hspace{0,2cm} -
	\dfrac{2(n-3)\mu^3b}{(n-1)(n-2)}\left\langle \left(E_2 \wedge_g E_i\right)(Z) , W\right\rangle , \\
 \left(C \cdot R\right)\left(E_i , E_j , Z , W; E_i , E_k\right)
 &=     0  , \\
 \left(C \cdot R\right)\left(E_i , E_j , Z , W; E_j , E_k\right)
 &=     0 .  \\
   \end{aligned}
	\right.
    \end{equation}

Thirdly, we compute the local components of the tensor
$ R \cdot C -  C\cdot R    $ on a  Wintgen   ideal submanifold.
For  any    index
   $ i    \in \left\{3 , \ldots ,  n\right\}$,	
\begin{equation}  \label{r.c-c.r_01}
\left\{
\begin{aligned}
 &\left(R \cdot C - C \cdot R\right)\left(E_1 , E_2 , Z , W ~; E_1 , E_i\right)
 =
 \dfrac{2(n-3)\mu^3a}{n-1}\left<\left(E_1 \wedge_g E_i\right)(Z) , W\right>\\
&  
  -
 \dfrac{2(n-3)\mu^2}{(n-1)(n-2)}
\left(\left(n-1\right)\left(H^2  + {\widetilde k}\right) + (n-2)\mu b + 2\mu^2\right)
                                  \left<\left(E_2 \wedge_g E_i\right)(Z) , W\right> , \\
&\left(R \cdot C - C \cdot R\right)\left(E_1 , E_2 , Z , W ~; E_2 , E_i\right)
 =
 -
 \dfrac{2n(n-3)\mu^3a}{(n-1)(n-2)}\left<\left(E_2 \wedge_g E_i\right)(Z) , W\right> \\ 
& 
+
 \dfrac{2(n-3)\mu^2}{(n-1)(n-2)}\left((n-1)\left(H^2
+ {\widetilde k}\right) - n\mu b 
- 2\mu^2\right)\left<\left(E_1 \wedge_g E_i\right)(Z) , W\right> .\\
 \end{aligned}
\right.
\end{equation}

\noindent       For  any   different indexes
   $ i $,  $ j     \in \left\{3 , \ldots ,  n\right\}$,	
\begin{equation}  \label{r.c-c.r_02}
\left\{
\begin{aligned}
&
\left(R\cdot C - C\cdot R\right)\left(E_1 , E_i , Z , W ~; E_1 , E_2\right) \\
= & \dfrac{4(n-3)\mu^3a}{n-1}\left\langle \left(E_1 \wedge_g E_i\right)(Z) , 
W\right\rangle
-
  \dfrac{4(n-3)\mu^3b}{n-1}\left\langle 
  \left(E_2 \wedge_g E_i\right)(Z) , W\right\rangle , \\
&
\left(R \cdot C - C \cdot R\right)\left(E_1 , E_i , Z , W ~; E_1 , E_j\right) \\
= &   \dfrac{2\mu^2}{(n-1)(n-2)}
\left((n-1)\left(H^2 + {\widetilde k}\right) + 2\mu b\right)
                               \left<\left(E_i \wedge_g E_j\right)(Z) , W\right> , \\
&
\left(R \cdot C - C \cdot R\right)\left(E_1 , E_i , Z , W ~; E_2 , E_i\right) \\
= & \dfrac{2(n-3)\mu^2}{(n-1)(n-2)}\left((n-1)\left(H^2
+ {\widetilde k}\right) -  n\mu b 
- 2\mu^2\right)\left<\left(E_1 \wedge_g E_2\right)(Z) , W\right> , \\
& \left(R \cdot C - C \cdot R\right)\left(E_1 , E_i , Z , W ~; E_i , E_j\right)
 = 0   . \\
 \end{aligned}
\right.
\end{equation}

\noindent       For  any   different indexes
   $ i $,  $ j     \in \left\{3 , \ldots ,  n\right\}$,		
\begin{equation}  \label{r.c-c.r_03}
\left\{
\begin{aligned}
& \left(R\cdot C - C\cdot R\right)\left(E_2 , E_i , Z , W ~; E_2 , E_1\right) \\
= &   \dfrac{4(n-3)\mu^3b}{n-1}\left\langle 
 \left(E_1 \wedge_g E_i\right)(Z) , W\right\rangle
+ \dfrac{4(n-3)\mu^3a}{n-1}\left\langle \left(E_2 \wedge_g E_i\right)(Z) , 
W\right\rangle , \\
& \left(R \cdot C - C \cdot R\right)\left(E_2 , E_i , Z , W ~; E_2 , E_j\right) \\
= &  \dfrac{2\mu^2}{(n-1)(n-2)}\left((n-1)
\left(H^2 + {\widetilde k}\right) 
- 2\mu b\right)\left<\left(E_i \wedge_g E_j\right)(Z) , W\right> , \\
&
\left(R \cdot C - C \cdot R\right)\left(E_2 , E_i , Z , W ~; E_i , E_1\right) \\
= & \dfrac{2(n-3)\mu^2}{(n-1)(n-2)}
\left((n-1)(H^2 + {\widetilde k}) + n\mu b\right)
               \left<\left(E_1 \wedge_g E_2\right)(Z) , W\right> , \\
& \left(R \cdot C - C \cdot R\right)\left(E_2 , E_i , Z , W ~; E_i , E_j\right)
 = 0 .  \\
 \end{aligned}
\right.
\end{equation}

\noindent       For  any   different indexes
   $ i $, $j$,  $ k     \in \left\{3 , \ldots ,  n\right\}$,	
\begin{equation}  \label{r.c-c.r_04}
\left\{
\begin{aligned}
\left(R \cdot C - C \cdot R\right)\left(E_i , E_j , Z , W ~; E_i , E_1\right)
 &=    \dfrac{2\mu^2}{(n-1)(n-2)}\left((n-1)\left(H^2 + {\widetilde k}\right) 
 + 2\mu b\right)
\left<\left(E_1 \wedge_g E_i\right)(Z) , W\right> \\
&\hspace{0,2cm}
+   \dfrac{4\mu^3a}{(n-1)(n-2)}\left<\left(E_2 \wedge_g E_i\right)(Z) , W\right> , \\
\left(R \cdot C - C \cdot R\right)\left(E_i , E_j , Z , W ~; E_i , E_2\right)
 &=  \dfrac{4\mu^3a}{(n-1)(n-2)}\left<\left(E_1 \wedge_g E_i\right)(Z) , W\right>
 \\
&\hspace{0,2cm}
+ \dfrac{2\mu^2}{(n-1)(n-2)}
\left((n-1)\left(H^2+{\widetilde k}\right)-2\mu b\right)
                  \left<\left(E_2 \wedge_g E_i\right)(Z), W\right> , \\
\left(R \cdot C - C \cdot R\right)\left(E_i , E_j , Z , W ~; E_i , E_k\right)
&=    0 , \\
\left(R \cdot C - C \cdot R\right)\left(E_i , E_j , Z , W ~; E_j , E_k\right)
 &=  0 .  \\
 \end{aligned}
\right.
\end{equation}



\noindent In terms of the equalities
  (\ref{r.c_01}), (\ref{r.c_02}), (\ref{r.c_03}), (\ref{r.c_04}),
   (\ref{c.r_01}), (\ref{c.r_02}), (\ref{c.r_03}), (\ref{c.r_04}),
and
(\ref{r.c-c.r_01}), (\ref{r.c-c.r_02}),
(\ref{r.c-c.r_03}), (\ref{r.c-c.r_04}), 
we prove theorems  1, 2, 3, 4 
and corollary  1 above.




\subsubsection{Proofs of theorems  {\sl \bf 5, 6, 7, 8, 9} 
and corollaries    {\sl \bf 2, 3}}

\noindent We compute the local components of the tensor $ Q(g , R) $ 
of a Wintgen ideal submanifold. Let   
be  tangent vector fields of $M$.
For   any index
   $ i    \in \left\{3 , \ldots ,  n\right\}$,
\begin{equation}\label{Qgr_01}
\left\{
\begin{aligned}
Q\left(g , R\right)\left(E_1 , E_2 , Z , W; E_1 , E_i\right)
 &= \left(b\mu   
 - 2\mu^2\right)\left\langle \left(E_2 \wedge_g E_i\right)(Z) , W\right\rangle ,  \\
Q\left(g , R\right)\left(E_1 , E_2 , Z , W; E_2 , E_i\right)
 &= \left(b\mu 
 +  2\mu^2\right)\left\langle \left(E_1 \wedge_g E_i\right)(Z) , W\right\rangle
 - a\mu\left\langle\left(E_2 \wedge_g E_i\right)(Z) , W\right\rangle  . \\
\end{aligned}
\right.
  \end{equation}

    \noindent   For   any    indexes
   $ i  $, $ j    \in \left\{3 , \ldots ,  n\right\}$	
	such that $ i \not= j$,
\begin{equation}\label{Qgr_02}
\left\{
\begin{aligned}
Q\left(g , R\right)\left(E_1 , E_i , Z , W; E_1 , E_2\right)
 &=    a\mu\left\langle \left(E_1 \wedge_g E_i\right)(Z) , W\right\rangle
  - 2 b\mu\left\langle \left(E_2 \wedge_g E_i\right)(Z) , W\right\rangle  , \\
Q\left(g , R\right)\left(E_1 , E_i , Z , W; E_1 , E_j\right)
 &=   b\mu\left\langle \left(E_i \wedge_g E_j\right)(Z) , W\right\rangle   ,  \\
Q\left(g , R\right)\left(E_1 , E_i , Z , W; E_2 , E_i\right)
 &=   \left(b\mu 
 + 2\mu^2\right)\left\langle \left(E_1 \wedge_g E_2\right)(Z) , W\right\rangle ,  \\
Q\left(g , R\right)\left(E_1 , E_i , Z , W; E_i , E_j\right)
  &=  0 . \\
 \end{aligned}
\right.
  \end{equation}

    \noindent   For   any indexes
   $ i  $, $ j    \in \left\{3 , \ldots ,  n\right\}$
		such that $ i \not= j$,
\begin{equation}\label{Qgr_03}
\left\{
\begin{aligned}
Q\left(g , R\right)\left(E_2 , E_i , Z , W; E_1 , E_2\right)
 &=    2b\mu \left\langle \left(E_1 \wedge_g E_i\right)(Z) , W\right\rangle
  - a\mu\left\langle \left(E_2 \wedge_g E_i\right)(Z) , W\right\rangle ,  \\
Q\left(g , R\right)\left(E_2 , E_i , Z , W; E_2 , E_j\right)
 &=    - b\mu \left\langle \left(E_i \wedge_g E_j\right)(Z) , W\right\rangle  , \\
Q\left(g , R\right)\left(E_2 , E_i , Z , W; E_i , E_1\right)
 &=   \left(-b\mu 
 + 2\mu^2\right)\left\langle \left(E_1 \wedge_g E_2\right)(Z) , W\right\rangle ,  \\
Q\left(g , R\right)\left(E_2 , E_i , Z , W; E_i , E_j\right)
 &=    0   .  \\
  \end{aligned}
\right.
  \end{equation}

    \noindent   For   any indexes
   $ i  $, $j$,  $ k    \in \left\{3 , \ldots ,  n\right\}$
	such that $ i \not= j$, $ i \not= k$ and $ j \not= k$,
\begin{equation}\label{Qgr_04}
\left\{
\begin{aligned}
Q\left(g , R\right)\left(E_i , E_j , Z , W; E_i , E_1\right)
   &=   b\mu\left\langle \left(E_1 \wedge_g E_j\right)(Z) , W\right\rangle , \\
Q\left(g , R\right)\left(E_i , E_j , Z , W; E_i , E_2\right)
   &=   a\mu\left\langle \left(E_1 \wedge_g E_j\right)(Z) , W\right\rangle
  -   b\mu\left\langle \left(E_2 \wedge_g E_j\right)(Z) , W\right\rangle  ,  \\
Q\left(g , R\right)\left(E_i , E_j , Z , W; E_i , E_k\right)
  &=  0 , \\
Q\left(g , R\right)\left(E_i , E_j , Z , W; E_j , E_k\right)
  &=  0 . \\
\end{aligned}
\right.
  \end{equation}

In terms of the equalities
  (\ref{r.c_01}), (\ref{r.c_02}), (\ref{r.c_03}), (\ref{r.c_04}),
  (\ref{c.r_01}), (\ref{c.r_02}), (\ref{c.r_03}), (\ref{c.r_04}),
(\ref{r.c-c.r_01}), (\ref{r.c-c.r_02}), (\ref{r.c-c.r_03}) , (\ref{r.c-c.r_04})
and
(\ref{Qgr_01}), (\ref{Qgr_02}), (\ref{Qgr_03}),   (\ref{Qgr_04}), we prove 
theorems  5, 6, 7, 8, 9
and corollaries  2, 3 above.

\subsubsection{Proofs of theorems   {\sl \bf 10, 11, 12}}

\noindent Now we compute the local components of the tensor $ Q(g , C) $
    for the considered Wintgen ideal
     submanifold. Let    $ X $, $ Y $,
 $ Z $, $ W $ be   tangent vector fields.
For   any index
   $ i   \in \left\{3 , \ldots ,  n\right\}$,
\begin{equation}\label{QgC_01}
\left\{
\begin{aligned}
  Q\left(g , C\right)\left(E_1 , E_2 , Z , W; E_1 , E_i\right)
&= - \dfrac{2(n-3)\mu^2}{n-2}\left\langle 
\left(E_2 \wedge_g E_i\right)(Z), W\right\rangle  , \\
  Q\left(g , C\right)\left(E_1 , E_2 , Z , W; E_2 , E_i\right)
 &=   \dfrac{2(n-3)\mu^2}{n-2}
 \left\langle \left(E_1 \wedge_g E_i\right)(Z) , W\right\rangle  .  \\
\end{aligned}
\right.
     \end{equation}

    \noindent   For   any indexes
   $ i  $, $ j    \in \left\{3 , \ldots ,  n\right\}$  such that $ i \not= j$,
\begin{equation}\label{QgC_02}
\left\{
\begin{aligned}
Q\left(g , C\right)\left(E_1 , E_i , Z , W; E_1 , E_2\right)
 &=   0   , \\
Q\left(g , C\right)\left(E_1 , E_i , Z , W; E_1 , E_j\right)
 &=       \dfrac{2\mu^2}{n-2}
 \left\langle \left(E_i \wedge_g E_j\right)(Z) , W\right\rangle  ,  \\
Q\left(g , C\right)\left(E_1 , E_i , Z , W; E_2 , E_i\right)
 &=   \dfrac{2(n-3)\mu^2}{n-2}
 \left\langle \left(E_1 \wedge_g E_2\right)(Z) , W\right\rangle  , \\
Q\left(g , C\right)\left(E_1 , E_i , Z , W; E_i , E_j\right)
 &=   0   . \\
   \end{aligned}
\right.
     \end{equation}

    \noindent   For   any indexes
   $ i  $, $ j    \in \left\{3 , \ldots ,  n\right\}$  such that $ i \not= j$,
\begin{equation}\label{QgC_03}
\left\{
\begin{aligned}
Q\left(g , C\right)\left(E_2 , E_i , Z , W; E_2 , E_1\right)
  &=  0 , \\
Q\left(g , C\right)\left(E_2 , E_i , Z , W; E_2 , E_j\right)
 &=   \dfrac{2\mu^2}{n-2}\left\langle 
 \left(E_i \wedge_g E_j\right)(Z) , W\right\rangle  ,  \\
Q\left(g , C\right)\left(E_2 , E_i , Z , W; E_i , E_1\right)
 &=      \dfrac{2(n-3)\mu^2}{n-2}
 \left\langle \left(E_1 \wedge_g E_2\right)(Z) , W\right\rangle  , \\
Q\left(g , C\right)\left(E_2 , E_i , Z , W; E_i , E_1\right)
 &=      0   .  \\
   \end{aligned}
\right.
     \end{equation}

    \noindent   For   any indexes
   $ i  $, $j$,  $ k    \in \left\{3 , \ldots ,  n\right\}$  such that $ i \not= j$,
$ i \not= k$ and  $ j \not= k$,
\begin{equation}\label{QgC_04}
\left\{
\begin{aligned}
Q\left(g , C\right)\left(E_i , E_j , Z , W; E_i , E_1\right)
 &=     \dfrac{2\mu^2}{n-2}\left\langle 
 \left(E_1 \wedge_g E_j\right)(Z) , W\right\rangle , \\
Q\left(g , C\right)\left(E_i , E_j , Z , W; E_i , E_2\right)
   &=    \dfrac{2\mu^2}{n-2}\left\langle 
   \left(E_2 \wedge_g E_j\right)(Z) , W\right\rangle ,  \\
Q\left(g , C\right)\left(X , Y , Z , W; E_i , E_k\right)
  &=  0  ,  \\
Q\left(g , C\right)\left(X , Y , Z , W; E_j , E_k\right)
  &=  0   . \\
   \end{aligned}
\right.
     \end{equation}

In terms of the equalities
  (\ref{r.c_01}), (\ref{r.c_02}), (\ref{r.c_03}), (\ref{r.c_04}),
  (\ref{c.r_01}), (\ref{c.r_02}), (\ref{c.r_03}), (\ref{c.r_04}),
(\ref{r.c-c.r_01}), (\ref{r.c-c.r_02}), (\ref{r.c-c.r_03}) , (\ref{r.c-c.r_04})
and
(\ref{QgC_01}),  (\ref{QgC_02}),  (\ref{QgC_03}),  (\ref{QgC_04}),
we prove theorems 10, 11, 12.

\subsubsection{Proofs of theorems   {\sl \bf 13, 14, 15, 16}
and corollary   {\sl \bf 4}}

\noindent We compute the local components of the tensor 
  $ Q(g ,   g\wedge\mbox{\rm Ricc}) $ of a Wintgen ideal submanifold.
 Let $X$, $Y$, $Z$, $W$ be
 tangent vector fields of $M$.
For any index
   $ i    \in \left\{3 , \ldots ,  n\right\}$,
  \begin{equation} \label{QgGvS01}
	\left\{
\begin{aligned}
Q(g ,  g\wedge\mbox{\rm Ricc})\left(E_1 , E_2 , Z , W ~; E_1 , E_i\right)
 &=
  \left((n-2)b\mu - 2\mu^2\right)\left<\left(E_2 \wedge_g E_i\right)(Z) , W\right> ,\\
Q(g ,  g\wedge\mbox{\rm Ricc})\left(E_1 , E_2 , Z , W ~; E_2 , E_i\right)
 &=
  \left((n-2)b\mu + 2\mu^2\right)\left<\left(E_1 \wedge_g E_i\right)(Z) , W\right> .\\
     \end{aligned}
    \right.
        \end{equation}

\noindent For any    different indexes
   $ i $, $ j   \in \left\{3 , \ldots ,  n\right\}$,
  \begin{equation} \label{QgGvS02}
	\left\{
\begin{aligned}
Q(g ,  g\wedge\mbox{\rm Ricc})\left(E_1 , E_i , Z , W ~; E_1 , E_2\right)
 &=
 - 2(n-2)b\mu\left<\left(E_2 \wedge_g E_i\right)(Z) , W\right> ,\\
Q(g ,  g\wedge\mbox{\rm Ricc})\left(E_1 , E_i , Z , W ~; E_1 , E_j\right)
 &= \left(- (n-2)b\mu + 2\mu^2\right)\left<\left(E_i\wedge_g E_j\right)(Z) , W\right> ,\\
Q(g ,  g\wedge\mbox{\rm Ricc})\left(E_1 , E_i , Z , W ~; E_2 , E_i\right)
&= \left((n-2)b\mu  +   2\mu^2\right)\left<\left(E_1\wedge_g E_2\right)(Z) , W\right> ,\\
Q(g ,  g\wedge\mbox{\rm Ricc})\left(E_1 , E_i , Z , W ~; E_i , E_j\right)
 &=    0 .\\
    \end{aligned}
    \right.
        \end{equation}

\noindent For any    different indexes
   $ i $, $ j   \in \left\{3 , \ldots ,  n\right\}$,
  \begin{equation} \label{QgGvS03}
	\left\{
\begin{aligned}
Q(g ,  g\wedge\mbox{\rm Ricc})\left(E_2 , E_i , Z , W ~; E_2 , E_1\right)
 &= 2(n-2)b\mu\left<\left(E_1 \wedge_g E_i\right)(Z) , W\right> ,\\
Q(g ,  g\wedge\mbox{\rm Ricc})\left(E_2 , E_i , Z , W ~; E_2 , E_j\right)
&=  \left(-(n-2)b\mu - 2\mu^2\right)\left<\left(E_i \wedge_g E_j\right)(Z) , W\right> ,\\
Q(g ,  g\wedge\mbox{\rm Ricc})\left(E_2 , E_i , Z , W ~; E_i , E_1\right)
&= \left(- (n-2)b\mu + 2\mu^2 \right)\left<\left(E_1 \wedge_g E_2\right)(Z) , W\right> ,\\
Q(g ,  g\wedge\mbox{\rm Ricc})\left(E_2 , E_i , Z , W ~; E_i , E_j\right)
 &=   0 .\\
    \end{aligned}
    \right.
        \end{equation}

\noindent
  For any    different indexes
   $ i $, $j$,  $ k   \in \left\{3 , \ldots ,  n\right\}$,
  \begin{equation} \label{QgGvS04}
	\left\{
\begin{aligned}
Q(g ,  g\wedge\mbox{\rm Ricc})\left(E_i , E_j , Z , W ~; E_i , E_1\right)
&=  \left((n-2)b\mu - 2\mu^2\right)\left<\left(E_1 \wedge_g E_j\right)(Z) , W\right> ,\\
Q(g ,  g\wedge\mbox{\rm Ricc})\left(E_i , E_j , Z , W ~; E_i , E_2\right)
&=  \left(- (n-2)b\mu  - 2\mu^2\right)\left<\left(E_2 \wedge_g E_i\right)(Z) , W\right> ,\\
Q(g ,  g\wedge\mbox{\rm Ricc})\left(E_i , E_j , Z , W ~; E_i , E_j\right)
 &=   	0 ,\\
Q(g ,  g\wedge\mbox{\rm Ricc})\left(E_i , E_j , Z , W ~; E_j , E_k\right)
 &=   	0 .\\
   \end{aligned}
    \right.
        \end{equation}

In terms of the equalities
  (\ref{r.c_01}), (\ref{r.c_02}), (\ref{r.c_03}), (\ref{r.c_04}),
  (\ref{c.r_01}), (\ref{c.r_02}), (\ref{c.r_03}), (\ref{c.r_04}),
(\ref{r.c-c.r_01}), (\ref{r.c-c.r_02}), (\ref{r.c-c.r_03}) , (\ref{r.c-c.r_04})
and
(\ref{QgGvS01}),  (\ref{QgGvS02}), (\ref{QgGvS03}), (\ref{QgGvS04}), we
proved theorems   13, 14, 15, 16 and corollary   4  above.

\subsubsection{Proofs of theorems   {\sl \bf 17, 18, 19, 20}
and corollary   {\sl \bf 5}}


  \noindent We compute the local components of the tensor $ Q(\mbox{\rm Ricc} ,   R) $
    for a Wintgen ideal submanifold.
Let $X$, $Y$, $Z$, $W$ be the
 tangent vector fields of $M$.
For   any index
   $ i     \in \left\{3 , \ldots ,  n\right\}$,
\begin{equation}\label{QSR_01}
\left\{
\begin{aligned}
&
Q\left(\mbox{\rm Ricc} , R\right)\left(E_1 , E_2 , Z , W; E_1 , E_i\right)
 =      a\mu\left(- \left(H^2 + {\widetilde k}\right) + 2\mu^2\right)
                  \left\langle \left(E_1 \wedge_g E_i\right)(Z) , W\right\rangle  \\
&
+ \left((n-2)(a^2 + b^2)\mu^2 - 2b\mu^3
- \left(2n - 4\right)\mu^2\left(H^2 + {\widetilde k}\right)
  +   b\mu\left(H^2 + {\widetilde k}\right)\right)
			      \left\langle \left(E_2 \wedge_g E_i\right)(Z) , W\right\rangle  ,   \\
&
Q\left(\mbox{\rm Ricc} , R\right)\left(E_1 , E_2 , Z , W; E_2 , E_i\right)
 =
  \left[- 2b\mu^3  - (n-2)(a^2 + b^2)\mu^2  \right.   \\
& \left.
 + (2n-4)\mu^2\left(H^2 + {\widetilde k}\right)
 +  b\mu\left(H^2 + {\widetilde k}\right)\right]
	                 \left\langle \left(E_1 \wedge_g E_i\right)(Z) , W\right\rangle  \\
&
+ \left(- 2a\mu^3 + a\mu\left(H^2 + {\widetilde k}\right)\right)
\left\langle \left(E_2 \wedge_g E_i\right)(Z) , W\right\rangle  .   \\						
   \end{aligned}
	\right.
\end{equation}

 \noindent  For   any different indexes
   $ i  $, $ j    \in \left\{3 , \ldots ,  n\right\}$,
\begin{equation}\label{QSR_02}
\left\{
\begin{aligned}
Q\left(\mbox{\rm Ricc} , R\right)\left(E_1 , E_i , Z , W; E_1 , E_2\right)
 &=   2a\mu \left(\left(H^2 + {\widetilde k}\right)  - 2\mu^2\right)
      \left\langle \left(E_1 \wedge_g E_i\right)(Z) , W\right\rangle  \\
&\hspace{0,2cm}
  + 2b\mu\left(- \left(H^2 + {\widetilde k}\right)  + 2\mu^2\right)
	          \left\langle \left(E_2 \wedge_g E_i\right)(Z) , W\right\rangle  ,   \\
Q\left(\mbox{\rm Ricc} , R\right)\left(E_1 , E_i , Z , W; E_1 , E_j\right)
 &=  \left(H^2 + {\widetilde k}\right)  \left(b\mu  +  2\mu^2\right)
           \left\langle \left(E_i \wedge_g E_j\right)(Z) , W\right\rangle   ,   \\
Q\left(\mbox{\rm Ricc} , R\right)\left(E_1 , E_i , Z , W; E_2 , E_i\right)
&= \left[- (n-2)\left(a^2 + b^2\right)\mu^2 - 2b\mu^3  \right.  \\
&\hspace{0,2cm}  \left.
+ (2n-4)\mu^2 (H^2 + {\widetilde k}) + b\mu(H^2 + {\widetilde k})\right]
           \left\langle \left(E_1 \wedge_g E_2\right)(Z) , W\right\rangle   ,   \\
Q\left(\mbox{\rm Ricc} , R\right)\left(E_1 , E_i , Z , W; E_i , E_j\right)
 &=   0 . \\
    \end{aligned}
	\right.
\end{equation}

 \noindent  For   any different  indexes
   $ i  $, $ j    \in \left\{3 , \ldots ,  n\right\}$,
\begin{equation}\label{QSR_03}
\left\{
\begin{aligned}
Q\left(\mbox{\rm Ricc} , R\right)\left(E_2 , E_i , Z , W; E_2 , E_1\right)
 &=     \left(2b\mu\left(H^2 + {\widetilde k}\right)  - 4b\mu^3\right)
      \left\langle \left(E_1 \wedge_g E_i\right)(Z) , W\right\rangle   \\
&\hspace{0,2cm}
  + \left(2a\mu\left(H^2 + {\widetilde k}\right)  - 4a\mu^3\right)
	          \left\langle \left(E_2 \wedge_g E_i\right)(Z) , W\right\rangle  ,   \\
Q\left(\mbox{\rm Ricc} , R\right)\left(E_2 , E_i , Z , W; E_2 , E_j\right)
&= \left(- b\mu  + 2\mu^2\right)\left(H^2 + {\widetilde k}\right)
 \left\langle \left(E_i \wedge_g E_j\right)(Z) , W\right\rangle ,\\
Q\left(\mbox{\rm Ricc} , R\right)\left(E_2 , E_i , Z , W; E_i , E_1\right)
  &=     \left[2b\mu^3	- (n-2)(a^2 + b^2)\mu^2  \right.   \\
&\hspace{0,2cm}  \left.
  +   (2n - 4)\mu^2\left(H^2 + {\widetilde k}\right)
		-    b\mu\left(H^2 + {\widetilde k}\right)\right]
 \left\langle \left(E_1 \wedge_g E_2\right)(Z) , W\right\rangle , \\
Q\left(\mbox{\rm Ricc} , R\right)\left(E_2 , E_i , Z , W; E_i , E_j\right)
  &=   0  . \\
   \end{aligned}
	\right.
\end{equation}

    \noindent   For  any different indexes
   $ i  $, $j$,  $ k    \in \left\{3 , \ldots ,  n\right\}$,
\begin{equation}\label{QSR_04}
\left\{
\begin{aligned}
Q\left(\mbox{\rm Ricc} , R\right)\left(E_i , E_j , Z , W; E_i , E_1\right)
  &=   \left(b\mu  + 2\mu^2\right)\left(H^2 + {\widetilde k}\right)
 \left\langle \left(E_1 \wedge_g E_j\right)(Z) , W\right\rangle  \\
&\hspace{0,2cm}
 +
  a\mu\left(H^2 + {\widetilde k}\right)\left\langle 
  \left(E_2 \wedge_g E_j\right)(Z) , W\right\rangle , \\
Q\left(\mbox{\rm Ricc} , R\right)\left(E_i , E_j , Z , W; E_i , E_2\right)
 &=   a\mu\left(H^2 + {\widetilde k}\right)
	 \left\langle \left(E_1 \wedge_g E_j\right)(Z) , W\right\rangle \\
&\hspace{0,2cm}
+ \left(- b\mu  + 2\mu^2\right)\left(H^2 + {\widetilde k}\right)
 \left\langle \left(E_2 \wedge_g E_j\right)(Z) , W\right\rangle ,\\
Q\left(\mbox{\rm Ricc} , R\right)\left(E_i , E_j , Z , W; E_i , E_k\right)
 &=   0 , \\
Q\left(\mbox{\rm Ricc} , R\right)\left(E_i , E_j , Z , W; E_j , E_k\right)
 &=   0 . \\
   \end{aligned}
	\right.
\end{equation}


From the equalities
  (\ref{r.c_01}), (\ref{r.c_02}), (\ref{r.c_03}), (\ref{r.c_04}),
  (\ref{c.r_01}), (\ref{c.r_02}), (\ref{c.r_03}), (\ref{c.r_04}),
(\ref{r.c-c.r_01}), (\ref{r.c-c.r_02}), (\ref{r.c-c.r_03}) , (\ref{r.c-c.r_04})
and
(\ref{QSR_01}), (\ref{QSR_02}), (\ref{QSR_03}), (\ref{QSR_04}),
we proved theorems   17, 18, 19, 20 and corollary 5 above.

\subsubsection{Proofs of theorems   {\sl \bf 21, 22, 23, 24}
and  corollary   {\sl \bf 6}}

  \noindent We compute the local components of the tensor $ Q(\mbox{\rm Ricc} ,   C) $
    for a Wintgen ideal submanifold.
Let $X$, $Y$, $Z$, $W$ be the
 tangent vector fields of $M$.
  For   any index
   $ i     \in \left\{3 , \ldots ,  n\right\}$,
\begin{equation}\label{QSC_01}
\left\{
\begin{aligned}
&
Q\left(\mbox{\rm Ricc} , C\right)\left(E_1 , E_2 , Z , W; E_1 , E_i\right) 
=   \dfrac{2(n-3)a\mu^3}{n-1}\left\langle 
 \left(E_1 \wedge_g E_i\right)(Z) , W\right\rangle  \\
&
 - \dfrac{2(n-3)\mu^2}{(n-1)(n-2)}  
\left((n-1)^2\left(H^2 + {\widetilde k}\right)  + (n -  2)b\mu  - 2\mu^2\right)  
 \left\langle  \left(E_2 \wedge_g E_i\right)(Z) , W\right\rangle    , \\
&
Q\left(\mbox{\rm Ricc} , C\right)\left(E_1 , E_2 , Z , W; E_2 , E_i\right) 
= 
-  \dfrac{2(n-3)a\mu^3}{n-1}\left\langle \left(E_2 \wedge_g E_i\right)(Z)  W\right\rangle     \\		
&
+ 
\dfrac{2(n-3)\mu^2}{(n-1)(n-2)}\left((n-1)^2\left(H^2 + {\widetilde k}\right)
 - (n -  2)b\mu  - 2\mu^2\right)   
\left\langle
 \left(E_1 \wedge_g E_i\right)(Z) , W\right\rangle .\\
   \end{aligned}
	\right.
\end{equation}

\noindent   For   any different  indexes
   $ i$, $ j  \in \left\{3 , \ldots ,  n\right\}$,
\begin{equation}\label{QSC_02}
\left\{
\begin{aligned}
&  Q\left(\mbox{\rm Ricc} , C\right)\left(E_1 , E_i , Z , W; E_1 , E_2\right) \\
=&   -  \dfrac{4(n-3)a\mu^3}{n-1}\left\langle 
 \left(E_1 \wedge_g E_i\right)(Z) , W\right\rangle
  + \dfrac{4(n-3)b\mu^3}{n-1}\left\langle 
  \left(E_2 \wedge_g E_i\right)(Z) , W\right\rangle ,    \\
& 
Q\left(\mbox{\rm Ricc} , C\right)\left(E_1 , E_i , Z , W; E_1 , E_j\right) \\
=&
  \dfrac{2\mu^2}{(n-1)(n-2)}\left((n-1)^2\left(H^2 + {\widetilde k}\right)
	+ 2(n -  2)b\mu  - 4\mu^2\right) 
   \left\langle \left(E_i \wedge_g E_j\right)(Z) , W\right\rangle    ,  \\
&
Q\left(\mbox{\rm Ricc} , C\right)\left(E_1 , E_i , Z , W; E_2 , E_i\right)  \\
=&    \dfrac{2(n-3)\mu^2}{(n-1)(n-2)}\left((n-1)^2\left(H^2 + {\widetilde k}\right)
 - (n -  2)b\mu  - 2\mu^2\right) 
 \left\langle 
 \left(E_1 \wedge_g E_2\right)(Z) , W\right\rangle ,\\
& Q\left(\mbox{\rm Ricc} , C\right)\left(E_1 , E_i , Z , W; E_i , E_j\right)  
=
 - \dfrac{2(n-3)a\mu^3}{n-1}
 \left\langle \left(E_1 \wedge_g E_2\right)(Z) , W\right\rangle   .  \\
   \end{aligned}
	\right.
\end{equation}

\noindent   For   any different  indexes
   $ i$, $ j  \in \left\{3 , \ldots ,  n\right\}$,
\begin{equation}\label{QSC_03}
\left\{
\begin{aligned}
&
Q\left(\mbox{\rm Ricc} , C\right)\left(E_2 , E_i , Z , W; E_2 , E_1\right) \\
= &    \dfrac{4(n-3)b\mu^3}{n-1}\left\langle 
 \left(E_1 \wedge_g E_i\right)(Z) , W\right\rangle
  - \dfrac{4(n-3)a\mu^3}{n-1}\left\langle 
  \left(E_2 \wedge_g E_i\right)(Z) , W\right\rangle   ,  \\
&
Q\left(\mbox{\rm Ricc} , C\right)\left(E_2 , E_i , Z , W; E_2 , E_j\right) \\
= &  \dfrac{2\mu^2}{(n-1)(n-2)}\left((n-1)^2\left(H^2 + {\widetilde k}\right)
 - 2(n -  2)b\mu  - 4\mu^2\right) 
\left\langle 
 \left(E_i \wedge_g E_j\right)(Z) , W\right\rangle    ,  \\
&
Q\left(\mbox{\rm Ricc} , C\right)\left(E_2 , E_i , Z , W; E_i , E_1\right) \\
= &    \dfrac{2(n-3)\mu^2}{(n-1)(n-2)}\left((n-1)^2\left(H^2 + {\widetilde k}\right)
 + (n -  2)b\mu  - 2\mu^2\right) 
 \left\langle 
 \left(E_1 \wedge_g E_2\right)(Z) , W\right\rangle , \\
&
Q\left(\mbox{\rm Ricc} , C\right)\left(E_2 , E_i , Z , W; E_i , E_j\right)
=    0   .\\
   \end{aligned}
	\right.
\end{equation}

\noindent   For   any different  indexes
   $ i$, $j$,  $ k  \in \left\{3 , \ldots ,  n\right\}$,
\begin{equation}\label{QSC_04}
\left\{
\begin{aligned}
&
Q\left(\mbox{\rm Ricc} , C\right)\left(E_i , E_j , Z , W; E_i , E_1\right)  
=  \dfrac{4a\mu^3}{n-1}\left\langle \left(E_2 \wedge_g E_j\right)(Z) , W\right\rangle \\
& +     \dfrac{2\mu^2}{(n-1)(n-2)}\left((n-1)^2\left(H^2 + {\widetilde k}\right)
 + 2(n -  2)b\mu  - 4\mu^2\right) 
 \left\langle \left(E_1 \wedge_g E_j\right)(Z) , W\right\rangle , \\
&
Q\left(\mbox{\rm Ricc} , C\right)\left(E_i , E_j , Z , W; E_i , E_2\right) 
= \dfrac{4a\mu^3}{n-1}\left\langle \left(E_1 \wedge_g E_j\right)(Z) , W\right\rangle \\
&
+ \dfrac{2\mu^2}{(n-1)(n-2)}\left((n-1)^2\left(H^2 + {\widetilde k}\right)
 - 2(n -  2)b\mu  - 4\mu^2\right) 
\left\langle \left(E_2 \wedge_g E_j\right)(Z) , W\right\rangle  ,\\
& Q\left(\mbox{\rm Ricc} , C\right)\left(E_i , E_j , Z , W; E_i , E_k\right)
 = 0  ,  \\
& Q\left(\mbox{\rm Ricc} , C\right)\left(E_i , E_j , Z , W; E_j , E_k\right)
 = 0   . \\
    \end{aligned}
	\right.
\end{equation}

From the equalities
  (\ref{r.c_01}), (\ref{r.c_02}), (\ref{r.c_03}), (\ref{r.c_04}),
  (\ref{c.r_01}), (\ref{c.r_02}), (\ref{c.r_03}), (\ref{c.r_04}),
(\ref{r.c-c.r_01}), (\ref{r.c-c.r_02}), (\ref{r.c-c.r_03}) , (\ref{r.c-c.r_04})
and
(\ref{QSC_01}), (\ref{QSC_02}), (\ref{QSC_03}), (\ref{QSC_04}),
we proved theorems    21, 22, 23, 24 and corollary  6 above.


\subsubsection{Proofs of theorems   {\sl \bf 25, 26, 27}
and corollaries {\sl \bf 7, 8, 9}}


  \noindent We compute the local components of the tensor
	$ Q(\mbox{\rm Ricc} ,   g\wedge\mbox{\rm Ricc}) $
    for a Wintgen ideal submanifold.
Let $X$, $Y$, $Z$, $W$ be the
 tangent vector fields of $M$.
For   any index
   $ i     \in \left\{3 , \ldots ,  n\right\}$,
\begin{equation}\label{QS_gS_01}
\left\{
\begin{aligned}
& 
Q\left(\mbox{\rm Ricc} , 
g\wedge\mbox{\rm Ricc}\right)\left(E_1 , E_2 , Z , W; E_1 , E_i\right) \\
= &  (n-1)a\mu\left[2(n-1)\left(H^2 + {\widetilde k}\right) + (n -  2)b\mu  - 2\mu^2\right]
 \left\langle\left(E_1 \wedge_g E_i\right)(Z) , W\right\rangle \\
&
+  \left[- 4\mu^4 + (n-2)^2b^2\mu^2   + 2(n-1)\mu^2\left(H^2 + {\widetilde k}\right)
- (n-1)(n -  2)b\mu  \left(H^2 + {\widetilde k}\right)\right]
  \left\langle \left(E_2 \wedge_g E_i\right)(Z) , W\right\rangle  ,\\
&Q\left(\mbox{\rm Ricc} , 
g\wedge\mbox{\rm Ricc}\right)\left(E_1 , E_2 , Z , W; E_2 , E_i\right)  \\
=&
\left[4\mu^4 -  (n -  2)^2b^2\mu^2 - 2(n-1)\mu^2
\left(H^2 + {\widetilde k}\right) \right. 
 \left.
- (n-1)(n-2)b\mu\left(H^2 + {\widetilde k}\right)\right]
\left\langle\left(E_1 \wedge_g E_i\right)(Z) , W\right\rangle  \\
& 
 -  (n-2)a\mu\left[2(n-1)
\left(H^2 + {\widetilde k}\right)
- (n -  2)b\mu  - 2\mu^2\right] 
\left\langle\left(E_2 \wedge_g E_i\right)(Z) , W\right\rangle . \\
     \end{aligned}
	\right.
\end{equation}

\noindent   For   any different indexes
   $ i $, $ j    \in \left\{3 , \ldots ,  n\right\}$,
\begin{equation}\label{QS_gS_02}
\left\{
\begin{aligned}
&
Q\left(\mbox{\rm Ricc} , g\wedge\mbox{\rm Ricc}\right)\left(E_1 , E_i , Z , W; E_1 , E_2\right)  \\
 = 
& -2(n-2)a\mu\left[2(n-1)\left(H^2 + {\widetilde k}\right) \right. 
\left.
 + (n -  2)b\mu - 2\mu^2\right]\left\langle \left(E_1 \wedge_g E_i\right)(Z) , W\right\rangle  \\
&
 +
2(n-1)(n-2)b\mu\left(H^2 + {\widetilde k}\right)
\left\langle\left(E_2 \wedge_g E_i\right)(Z) , W\right\rangle  ,\\
&
Q\left(\mbox{\rm Ricc} , g\wedge\mbox{\rm Ricc}\right)\left(E_1 , E_i , Z , W; E_1 , E_j\right) \\
= &
 \left[- (n-1)(n -  2)b\mu\left(H^2 + {\widetilde k}\right)  \right. 
\left.
+ 2(n-1)\mu^2\left(H^2 + {\widetilde k}\right)\right]
\left\langle \left(E_i \wedge_g E_j\right)(Z) , W\right\rangle , \\
&
Q\left(\mbox{\rm Ricc} , g\wedge\mbox{\rm Ricc}\right)\left(E_1 , E_i , Z , W; E_2 , E_j\right)  \\
 = &  \left[4\mu^4 - (n -  2)^2b^2\mu^2
 - 2(n-1)\mu^2\left(H^2 + {\widetilde k}\right)\right. 
\left.
- (n-1)(n -  2)b\mu\left(H^2 + {\widetilde k}\right)\right]
\left\langle \left(E_1 \wedge_g E_2\right)(Z) , W\right\rangle , \\
&Q\left(\mbox{\rm Ricc} , g\wedge\mbox{\rm Ricc}\right)\left(E_1 , E_i , Z , W; E_i , E_j\right)
 =  0  .  \\
     \end{aligned}
	\right.
\end{equation}

\noindent   For   any different  indexes
   $ i $, $ j    \in \left\{3 , \ldots ,  n\right\}$,
\begin{equation}\label{QS_gS_03}
\left\{
\begin{aligned}
 &
Q\left(\mbox{\rm Ricc} , g\wedge\mbox{\rm Ricc}\right)\left(E_2 , E_i , Z , W; E_2 , E_1\right) \\
= &
 - 2(n-1)(n-2)b\mu\left(H^2 + {\widetilde k}\right)
\left\langle \left(E_1 \wedge_g E_i\right)(Z) , W\right\rangle    \\
&
- 2(n-2)a\mu\left[2(n-1)\left(H^2 + {\widetilde k}\right)  - (n-2)b\mu - 2\mu^2\right]
\left\langle \left(E_2 \wedge_g E_i\right)(Z) , W\right\rangle  ,  \\
&
Q\left(\mbox{\rm Ricc} , g\wedge\mbox{\rm Ricc}\right)\left(E_2 , E_i , Z , W; E_2 , E_j\right) \\
= & \left[(n-1)(n-2)b\mu\left(H^2 + {\widetilde k}\right)
 \right.   
  \left.
+ 2(n-1)\mu^2\left(H^2 + {\widetilde k}\right)\right]\left\langle \left(E_i \wedge_g E_j\right)(Z) , W\right\rangle  ,  \\
&
Q\left(\mbox{\rm Ricc} , g\wedge\mbox{\rm Ricc}\right)\left(E_2 , E_i , Z , W; E_i , E_1\right) \\
= & \left[4\mu^4 - (n-2)b^2\mu^2 \right.  
 \left.
+ (n-1)(n-2)b\mu\left(H^2 + {\widetilde k}\right)
 -  2(n-1)\mu^2\left(H^2 + {\widetilde k}\right)\right]
\left\langle \left(E_1 \wedge_g E_2\right)(Z) , W\right\rangle    \\
&Q\left(\mbox{\rm Ricc} , g\wedge\mbox{\rm Ricc}\right)\left(E_2 , E_i , Z , W; E_i , E_j\right)
= 0  .   \\
     \end{aligned}
	\right.
\end{equation}

\noindent   For   any  different indexes
   $ i $, $j$, $ k    \in \left\{3 , \ldots ,  n\right\}$,
\begin{equation}\label{QS_gS_04}
\left\{
\begin{aligned}
 &Q\left(\mbox{\rm Ricc} , g\wedge\mbox{\rm Ricc}\right)\left(E_i , E_j , Z , W; E_i , E_1\right)  \\
= &
\left[(n-1)\left(H^2 + {\widetilde k}\right)\left(- (n-2)b\mu + 2\mu^2 \right)\right]
\left\langle \left(E_1 \wedge_g E_j\right)(Z) , W\right\rangle    \\
&
- 2(n-1)(n-2)a\mu\left(H^2 + {\widetilde k}\right)
\left\langle \left(E_2 \wedge_g E_j\right)(Z) , W\right\rangle   , \\
&
Q\left(\mbox{\rm Ricc} , g\wedge\mbox{\rm Ricc}\right)\left(E_i , E_j , Z , W; E_i , E_2\right) \\
=
& - 2(n-1)(n-2)a\mu\left(H^2 + {\widetilde k}\right)
\left\langle \left(E_1 \wedge_g E_j\right)(Z) , W\right\rangle    \\
&
+ \left[(n-1)(n+2)b\mu\left(H^2 + {\widetilde k}\right)
+ 2(n-1)\mu^2\left(H^2 + {\widetilde k}\right)\right]
\left\langle \left(E_2 \wedge_g E_j\right)(Z) , W\right\rangle  ,  \\
&
Q\left(\mbox{\rm Ricc} , g\wedge\mbox{\rm Ricc}\right)\left(E_i , E_j , Z , W; E_i , E_k\right) =  0  , \\
&
Q\left(\mbox{\rm Ricc} , g\wedge\mbox{\rm Ricc}\right)\left(E_i , E_j , Z , W; E_j , E_k\right) =  0 .  \\
     \end{aligned}
	\right.
\end{equation}

From the equalities
  (\ref{r.c_01}), (\ref{r.c_02}), (\ref{r.c_03}), (\ref{r.c_04}),
  (\ref{c.r_01}), (\ref{c.r_02}), (\ref{c.r_03}), (\ref{c.r_04}),
(\ref{r.c-c.r_01}), (\ref{r.c-c.r_02}), (\ref{r.c-c.r_03}) , (\ref{r.c-c.r_04})
and
(\ref{QS_gS_01}), (\ref{QS_gS_02}), (\ref{QS_gS_03}), (\ref{QS_gS_04}),
we proved theorems  25, 26, 27
and corollaries 7, 8, 9 above.

\vspace{5mm}

\noindent
{\bf{Funding.}} 
The first two authors of this paper are supported 
by the Wroc\l aw University of Environmental and Life Sciences (Poland).



\noindent
\footnotesize{Ryszard Deszcz\\
retired employee of the Department of Applied Mathematics\\
Wroc\l aw University of Environmental and Life Sciences\\
Grunwaldzka 53, 50-357 Wroc\l aw, Poland}\\
{\sf E-mail: Ryszard.Deszcz@upwr.edu.pl}\\
\textbf{ORCID ID: 0000-0002-5133-5455}
\newline

\noindent
\footnotesize{Ma\l gorzata G\l ogowska\\
Department of Applied Mathematics \\
Wroc\l aw University of Environmental and Life Sciences\\
Grunwaldzka 53, 50-357 Wroc\l aw, Poland}\\
{\sf E-mail: Malgorzata.Glogowska@upwr.edu.pl}\\
\textbf{ORCID ID: 0000-0002-4881-9141}
\newline

\noindent
\footnotesize{Miroslava Petrovi\'{c}-Torga\v{s}ev\\
Department of Sciences and Mathematics\\
State University of Novi Pazar\\
Vuka Karad\v{z}i\'{c}a 9\\
36300 Novi Pazar, Serbia}\\
{\sf E-mail: mirapt@kg.ac.rs}\\
\textbf{ORCID ID: 0000-0002-9140-833X}
\newline

\noindent
\footnotesize{Georges Zafindratafa\\
Professeur \'{E}m\'{e}rite of the
Laboratoire
de Math\'{e}matiques pour l'Ing\'{e}nieur (LMI)\\
Universit\'{e} Polytechnique Hauts-de-France\\
59313 Va\-len\-cien\-nes cedex 9, France}\\
{\sf E-mail: Georges.Zafindratafa@uphf.fr}\\
\textbf{ORCID ID: 0009-0001-7618-4606}


\begin{thebibliography}{99}


\bibitem{ECA} 
E. Cartan,
    \emph{Le{\c{c}}ons sur la g{\'e}om{\'e}trie des espaces de Riemann},
2nd ed., Paris: Gauthier-Villars, ($1946$).


\bibitem{BYCH1993} 
B.-Y. Chen,
  \emph{Some pinching and classification  theorems for minimal submanifolds},
Arch. f{\"{u}}r Math.   $\mathbf{60}$
 ($1993$), 568--578.
 


\bibitem{Ch11} 
B.-Y. Chen, 
  \emph{Complex extensors and Lagrangian submanifolds
	in complex Euclidean  spaces},
T\^{o}hoku Math. J.
$\mathbf{49}$ ($1997$),  277--297.


\bibitem{BYCH2008} 
B.-Y. Chen,
  \emph{$ \delta $-invariants, inequalities of
  submanifolds  and their applications},
  Topics in Differential Geo\-met\-ry, Ch. 2,
  Editors A. Mihai, I. Mihai and R. Miron, Editura Academiei Rom\^{a}ne,
  Bucharest, $2008$.


\bibitem{11} 
B.-Y. Chen, 
  \emph{Classification of Wintgen ideal surfaces
	in Euclidean $4$-space with equal Gauss and normal curvatures},
Ann. Glob Anal. Geom.
$\mathbf{38}$ ($2010$),  145--160.

\bibitem{Ch12} 
B.-Y. Chen, 
\emph{A Wintgen type inequality for surfaces in $4D$ neutral
	neutral pseudo-Riemannian space forms and its applications
	to minimal immersions},
Int. J. Math. Sci.
$\mathbf{1}$ ($2010$),  1--12.



\bibitem{BYCH2011} 
B.-Y. Chen,
  \emph{Pseudo-Riemannian geometry, $ \delta $-invariants and applications},
World Sci., Hackensack, New Jersey,
 ($2011$).



\bibitem{Ch13} 
B.-Y. Chen, 
  \emph{Wintgen ideal surfaces in four-dimensional neutral
	indefinite   space form $ \R^4_2(c) $},
Results Math.
$\mathbf{61}$ ($2012$),  329--345.


\bibitem{CHEN-2017}
B.-Y. Chen,
\emph{Differential Geometry of Warped Product Manifolds and Submanifolds},
World Sci., 2017.


\bibitem{BYCH2021} 
B.-Y. Chen,
\emph{Recent Developments in Wintgen Inequality and Wintgen Ideal submanifolds},
Int. Electron. J. Geom.
$\mathbf{14}$  ($2021$), 1--40.


\bibitem{3} 
T. Choi and Z. Lu,
\emph{On the DDVV conjecture and the comass in calibrated geometry}, I.
Math. Z.
$\mathbf{260}$ ($2008$), 409--429.


\bibitem{14} 
S. Decu, M. Petrovi\'{c}-Torga\v{s}ev,  A. \v{S}ebekovi\'{c} and  L. Verstraelen,
\emph{On the intrinsic Deszcz  symmetries
and the extrinsic Chen character of Wintgen  ideal submanifolds},
Tamkang J. Math.
$\mathbf{41}$ ($2010$), 109--116.


\bibitem{15} 
S. Decu, M. Petrovi\'{c}-Torga\v{s}ev,  A. \v{S}ebekovi\'{c} and  L. Verstraelen,
\emph{On the Roter type of Wintgen ideal submanifolds},
Rev. Roumaine Math. Pures Appl.
$\mathbf{57}$ ($2012$), 75--90.



\bibitem{1} 
P. J. De Smet, F. Dillen, L. Verstraelen  and L. Vrancken,
\emph{A pointwise inequality in submanifold theory},
Arch. Math. (Brno)
$\mathbf{35}$ ($1999$), 115--128.


\bibitem{DES2}
 R.  Deszcz,
 \emph{On pseudo-symmetric spaces}, Bull. Belg. Math. Soc., Ser. A,
 $\mathbf{44}$  ($1992$), 1--34.


\bibitem{DGHHY} 
R. Deszcz, M. G\l ogowska, H. Hashiguchi, M. Hotlo\'{s} and M. Yawata,
\emph{On semi-Riemannian manifolds satisfying some conformally invariant curvature condition}, 
Colloq. Math. {\bf{131}} (2013), 149--170.



\bibitem{2023_DGHP-TZ}
R. Deszcz, M. G\l ogowska, M. Hotlo\'{s}, M. Petrovi\'{c}-Torga\v{s}ev and G. Zafindratafa,
\emph{A note on some gene\-ra\-li\-zed curvature tensor},
Int. Electron. J. Geom. 16 (1) (2023), 379--397.


\bibitem{2023_b_DGHP-TZ}
R. Deszcz, M. G\l ogowska, M. Hotlo\'{s}, M. Petrovi\'{c}-Torga\v{s}ev and G. Zafindratafa,
\emph{On semi-Riemannian manifolds satisfying some generalized Einstein metric conditions},
Int. Electron. J. Geom. 16 (2) (2023), 539--576.


\bibitem{DGHS2011}
R. Deszcz, M. G\l ogowska, M. Hotlo\'{s}, and K. Sawicz,
\emph{Survey on Generalized Einstein Metric Conditions},
in: Advances in Lorentzian Geometry:
Proceedings of the Lorentzian Geometry Conference in Berlin,
AMS/IP Studies in Advanced Mathematics $\mathbf{49}$, S.-T. Yau (series ed.),
M. Plaue, A.D. Rendall and M. Scherfner (eds.), ($2011$), 27--46.


\bibitem{DGHS-2022}
R. Deszcz, M. G\l ogowska, M. Hotlo\'{s} and Z. \c Sent\"{u}rk,
\emph{On some quasi-Einstein and $2$-quasi-Einstein ma\-ni\-folds},
AIP Conference Proceedings 2483, 100001 (2022);
https://doi.org/10.1063/5.0118057


\bibitem{DGPV}
R. Deszcz,   M. G\l ogowska,  M. Petrovi\'{c}-Torga\v{s}ev and L. Verstraelen,
\emph{On the Roter  type of Chen  ideal submanifolds},
Results Math. {\bf 59} (2011), 401--413.


\bibitem{DGP-TV02}
R. Deszcz, M. G\l ogowska, M. Petrovi\'{c}-Torga\v{s}ev and L. Verstraelen,
\emph{Curvature properties of some class of minimal hypersurfaces in Euclidean spaces},
Filomat {\bf{29}} (2015), 479--492.


\bibitem{DGPSS}
R. Deszcz, M. G\l ogowska, M. Plaue, K. Sawicz and M. Scherfner,
\emph{On hypersurfaces in space forms satisfying particular curvature conditions of Tachibana type},
Kragujevac J. Math.
$\mathbf{35}$ ($2011$), 223--247.


\bibitem{DESZ}
 R.  Deszcz,  S. Haesen  and L. Verstraelen,
 \emph{On natural symmetries},
 Topics in Differential Geometry, Ch. 6,
  Editors A. Mihai, I. Mihai and R. Miron, Editura Academiei Rom\^{a}ne,
  Bucharest, $2008$.


 
\bibitem{2008_DP-TVZ} 
R. Deszcz, M. Petrovi\'{c}-Torga\v{s}ev, L. Verstraelen and G. Zafindratafa,
\emph{On the intrinsic symmetries of Chen ideal submanifolds},
Bull. Transilvania Univ. Brasov,
Ser. Math., Informatics, Physics, {\bf{15}} ({\bf{50}}) (2008), 99--108. 


\bibitem{DP-TVZ}
R. Deszcz, M. Petrovi\'{c}-Torga\v{s}ev, L. Verstraelen and G. Zafindratafa,
\emph{On Chen ideal submanifolds satisfying some conditions of pseudo-symmetry type},
Bull. Malaysian Math. Sci. Soc.
$\mathbf{39}$ ($2016$), 103--131.


\bibitem{13} 
R. Deszcz, M. Petrovi\'{c}-Torga\v{s}ev,  Z. {\c{S}}ent\"{u}rk and  L. Verstraelen,
\emph{Characterization of the  pseudo-symmetries of Wintgen ideal submanifolds of dimension} $3$,
Publ. Inst. Math. (Beograd) (N.S.)
$\mathbf{88}$ ($\mathbf{102}$) ($2010$), 53--65.


\bibitem{44} R. Deszcz, L. Verstraelen and \c S. Yaprak, 
\emph{Pseudosymmetric hypersurfaces in $4$-dimensional spaces of constant curvature}, 
Bull. Inst. Math. Acad. Sinica {\bf{22}} (1994), 167--179.


\bibitem{DeVerYap}
R. Deszcz, L. Verstraelen and \c S. Yaprak,
\emph{On 2-quasi-umbilical hypersurfaces in conformally flat spaces},
Acta Math. Hung. 78 (1998), 45--57.


\bibitem{DY} R. Deszcz and \c S. Yaprak,
\emph{Curvature properties of Cartan hypersurfaces}, 
Colloq. Math. {\bf{67}} (1994), 91--98.


\bibitem{4} 
J. Ge and Z. Tang,
\emph{A proof of the DDVV conjecture and its equality case},
Pacific J. Math.
$\mathbf{237}$ ($2008$), 87--95.


\bibitem{5} 
J. Ge and Z. Tang,
\emph{A survey on the DDVV conjecture}, arXiv: 1006.5326v1
 ($2010$).


\bibitem{GuaRod} 
I. V.  Gudalupe and I. Rodriguez,
\emph{Normal curvature of surfaces in space forms},
Pacific J. Math.
$\mathbf{106}$ ($1983$), 95--103.


\bibitem{HaVer03}
S.  Haesen  and  L. Verstraelen,
\emph{Classification of the pseudo-symmetric space-times},
J. Math. Phys.
$\mathbf{45}$ ($2004$),  2343--2346.


\bibitem{HaVer01}
S.  Haesen and L. Verstraelen,
\emph{Properties of a scalar curvature invariant depending on two planes},
Manuscripta Math.
$\mathbf{122}$ ($2007$), 59--72.


\bibitem{HaVer02} 
S.  Haesen  and  L.  Verstraelen,
\emph{Natural Intrinsic Geometrical Symmetries},
Symmetry, Integrability and Geometry: Methods and Applications SIGMA
$\mathbf{5}$ ($2009$), $086$, 15 pages.


\bibitem{2} 
Z. Lu,
\emph{Normal scalar curvature conjecture and its applications},
J. Funct. Anal.
$\mathbf{261}$ ($2011$), 1284--1308.


\bibitem{MPTAP}
M. Petrovi\'{c}-Torga\v{s}ev and A. Panti\'{c},
\emph{Pseudo-symmetries of generalized Wintgen ideal Lagrangian sub\-ma\-ni\-folds},
Publ. Inst. Math. (Beograd) (N.S.)
${\mathbf{103}}$ ({$\mathbf{117}$}) ($2018$), 181--190.



\bibitem{12} 
M. Petrovi\'{c}-Torga\v{s}ev and  L. Verstraelen,
\emph{On Deszcz symmetries of Wintgen ideal submanifolds},
Arch. Math. (Brno),
$\mathbf{44}$  ($2008$),  57--67.

\bibitem{ROUX} 
B. Rouxel, 
\emph{Sur une famille de A-surfaces d'un espace euclidien} $E^4$,
Proc. 10 \"{O}sterreichischer Mathematiker Kongress, Innsbruck, 1981, 185.


\bibitem{SEB}
A. \v{S}ebekovi\'{c},
\emph{Symmetries of Wintgen ideal submanifolds},
Bull. Transilv. Univ. Bra\c{s}ov,
Ser. Math., Informatics, Physics,
ser. III, $\mathbf{1}$ ($\mathbf{50}$)  ($2008$), 333--341.


\bibitem{SEBMPTAP}
A. \v{S}ebekovi\'{c}, M. Petrovi\'{c}-Torga\v{s}ev and A. Panti\'{c},
\emph{Pseudosymmetry properties of generalized Wintgen ideal Legendrian submanifolds},
Filomat $\mathbf{33}$ ($2019$), 1209--1215.


\bibitem{SEN}
Z. {\c{S}}ent\"{u}rk,
\emph{Characterisation of the Deszcz symmetric ideal Wintgen submanifolds},
An. Stiint Univ. Al. I. Cuza Ia\c si Mat. (N.S.) $\mathbf{53}$ ($2007$),
suppl. I, 309--316.


\bibitem{Sz1982} Z. I. Szab\'o,
 \emph{Structure theorems on Riemannian spaces satisfying}
    $R(X,Y)\cdot R = 0$.
    {\em I. The local version},
J. Differential Geom.
$\mathbf{17}$ ($1982$), 531--582.



\bibitem{LV1}
L. Verstraelen,
\emph{Comments on the pseudo-symmetry in the sense of Ryszard Deszcz},
in: Geometry and Topology of Submanifolds, $\mathbf{VI}$,
World Sci., Singapore, 1994,
119--209.


\bibitem{V2}
L. Verstraelen,
\emph{A coincise mini history of Geometry},
Kragujevac J. Math. 38 (2014), 5--21.


\bibitem{LV2}
L. Verstraelen,
\emph{Natural extrinsic geometrical symmetries - an introduction -},
in: Recent Advances in the Geometry of Submanifolds
Dedicated to the Memory of Franki Dillen (1963-2013),
in: Contemporary Mathematics, $\mathbf{674}$ (2016), pp. 5--16.


\bibitem{LV3-Foreword}
L. Verstraelen,
\emph{Foreword}, in:  B.-Y. Chen,
\emph{Differential Geometry of Warped Product Manifolds and Submanifolds},
World Scientific, 2017, vii--xxi.


\bibitem{LV4}
L. Verstraelen,
\emph{Submanifolds theory--A contemplation of submanifolds},
in: Geometry of Submanifolds,
AMS Special Session on Geometry of Submanifolds
in Honor of Bang-Yen Chen's 75th Birthday,
October 20-21, 2018,
University of Michigan, Ann Arbor, Michigan,
J. Van der Veken et al. (eds.),
Contemporary Math. 756, Amer. Math. Soc., 2020, 21--56.


\bibitem{Vilcu-2022}
G.-E. Vilcu, 
\emph{Curvature Inequalities for Slant Submanifolds in Pointwise Kenmotsu Space Forms},
in: 
Contact Geometry of Slant Submanifolds,
B.-Y. Chen, M.H. Shahid, F. Al-Solamy (eds.), Singapore: Springer. 13--37 (2022). 



\bibitem{8} 
P. Wintgen,
\emph{Sur l'in{\' e}galit{\' e} de Chen-Willmore},
C. R.  Acad. Sci. Paris
$\mathbf{288}$ ($1979$),  993--995.

\end{thebibliography}
     \end{document}